\documentclass[onefignum,onetabnum]{siamart171218}

\usepackage{lipsum}
\usepackage{amsfonts}
\usepackage{graphicx}
\usepackage{epstopdf}
\ifpdf
  \DeclareGraphicsExtensions{.eps,.pdf,.png,.jpg}
\else
  \DeclareGraphicsExtensions{.eps}
\fi

\newsiamremark{remark}{Remark}
\newsiamremark{hypothesis}{Hypothesis}
\crefname{hypothesis}{Hypothesis}{Hypotheses}
\newsiamthm{claim}{Claim}

\headers{RTSMS}{B. Hashemi and Y. Nakatsukasa }

\title{RTSMS: Randomized Tucker with single-mode sketching\thanks{Version of \today.
}}

\author{Behnam Hashemi\thanks{School of Computing and Mathematical Sciences, University of Leicester, Leicester, LE1 7RH, UK (\email{b.hashemi@le.ac.uk}).}
\and 
Yuji Nakatsukasa\thanks{Mathematical Institute, University of Oxford, Oxford, OX2 6GG, UK\\ (\email{nakatsukasa@maths.ox.ac.uk}).}
}

\usepackage{amsopn}

\DeclareMathOperator*{\minimize}{minimize}
\DeclareMathOperator*{\minimizebf}{\bf{minimize}}

\interfootnotelinepenalty=10000 \usepackage{tablefootnote}

 \newtheorem{example}{Example}
\usepackage{graphicx}
 \graphicspath{{figures/}} 

\usepackage[normalem]{ulem}
\usepackage{multirow}
\usepackage{epstopdf}
\usepackage{listings}
\usepackage[]{mcode}
\usepackage{algorithm}
\usepackage{latexsym}
\usepackage{showidx}
\usepackage{latexsym}
\usepackage{amssymb}
\renewcommand{\vec}{{\rm vec}}
\newcommand{\eqnref}[1]{(\ref{#1})}
\newcommand{\ignore}[1]{}

\usepackage{todonotes}

\usepackage{upquote}
\usepackage{graphicx,color}
\usepackage{algcompatible}
\usepackage{latexsym}
\usepackage{amssymb}
\usepackage{amsmath}
\usepackage{pdflscape}

\usepackage{mathrsfs}
\usepackage{soul}
\usepackage{bigfoot}
\usepackage{stmaryrd} \def\RTSMS{RTSMS}

\usepackage[normalem]{ulem}

\begin{document}
\maketitle

\begin{abstract}
We propose RTSMS (Randomized Tucker via Single-Mode-Sketching), a randomized algorithm for approximately computing a low-rank Tucker decomposition of a given tensor. It uses sketching and least-squares to compute the Tucker decomposition in a sequentially truncated manner. The algorithm only sketches one mode at a time, so the sketch matrices are significantly smaller than alternative approaches. The algorithm is demonstrated to be competitive with existing methods, sometimes outperforming them by a large margin.
\end{abstract}



\section{Introduction}
The Tucker decomposition is a family of representations that break up a given tensor $\mathcal{A} \in\mathbb{R}^{n_1\times n_2\times \cdots \times n_d}$ into the multilinear product of a core tensor $\mathcal{C} \in\mathbb{R}^{r_1\times r_2\times \cdots \times r_d}$ and a factor matrix $F_k \in\mathbb{R}^{n_k\times r_k}$ $(r_k\leq n_k)$
along each mode $k=1,2,\dots d$, i.e., 
\[
\mathcal{A} = \mathcal{C} \times_1 F_1 \times_2 F_2 \dots \times_d F_d,
\]
see Section~\ref{prelim:sec} for the definition of mode-$k$ product $\times_k$. Assuming that $\mathcal{A}$ can be well approximated by a low multilinear rank decomposition ($r_k\ll n_k$ for some or all $k$), one takes advantage of the fact that the core tensor $\mathcal{C}$ can be significantly smaller than $\mathcal{A}$. The canonical polyadic (CP) decomposition, also very popular in multilinear data analysis, is a special case of the Tucker decomposition in which the core tensor $\mathcal{C}$ has to be diagonal. 

The history goes back to the 1960s when Tucker introduced the concept as a tool in quantitative psychology \cite{Tucker63}, as well as algorithms \cite{Tucker66} for its computation. The decomposition has various applications such as dimensionality reduction, face recognition, image compression~
\cite{Pan21,VasilescuThesis}, etc.
Several deterministic and randomized algorithms have been developed for the computation of Tucker decomposition~\cite{Bucci23, Che20, de2000multilinear, Malik18, Minster23, Minster20, Sun20, Vannieuwenhoven12, Zhou14}. 

In what follows we give a very short outline of certain aspects of Tucker decomposition and refer the reader to the review \cite{kolda2009tensor} by Kolda and Bader for other aspects including citations to several contributions. The higher order orthogonal iteration (HOOI) by De Lathauwer, De Moor and Vandewalle~\cite{de2000best}
is an alternating least squares (ALS) method that uses the SVD to find the best multilinear rank ${\bf r}$ approximation of $\mathcal{A}$. The same authors \cite{de2000multilinear} introduced a generalization of the singular value decomposition of matrices to tensors called HOSVD. Vannieuwenhoven, Vandebril and Meerbergen \cite{Vannieuwenhoven12} introduced an efficient algorithm called STHOSVD which computes the HOSVD via a sequential truncation of the tensor taking advantage of the compressions achieved when processing all of the previous modes.

Cross algorithms constitute a wide class of powerful deterministic algorithms for computing a Tucker decomposition.\ These algorithms compute Tucker decompositions by interpolating the given tensor on a carefully selected set of pivot elements, also refereed to as ``cross'' points. Notable examples include the Cross-3D algorithm of Oseledets, Savostianov and Tyrtyshnikov \cite{Oseledets08cross3D} and the fiber sampling methods of Caiafa and Cichocki \cite{Caiafa10}. Relevant Tucker decomposition algorithms in the context of 3D function approximation include the slice-Tucker decomposition \cite{hashemichebfun3} and the fiber sampling algorithm of Dolgov, Kressner and Str\"ossner \cite{dolgov21} which relies on oblique projections using subindices chosen based on the discrete empirical interpolation (DEIM) \cite{chaturantabut2010nonlinear} method.

In this paper, we introduce RTSMS (Randomized Tucker via Single-Mode-Sketching), which is based on sequentially finding the factor matrices in the Tucker approximation while truncating the tensor. In this sense it is 
similar to STHOSVD~\cite{Vannieuwenhoven12}. 
Crucially, RTSMS is based on sketching the tensor from just one side (mode) at a time, not two (or more)---thus avoiding an operation that can be a computational bottleneck. 
This is achieved by finding a low-rank approximation of unfoldings based on generalized Nystr\"om (GN) (instead of randomized SVD~\cite{HMT}), but replacing the second 
sketch with a subsampling matrix chosen via the leverage scores of the first sketch. In addition, we employ regularization and iterative refinement techniques to improve the numerical stability.

Another key aspect of \RTSMS\ is its ability to find the rank adaptively, given a required error tolerance. We do this by blending matrix rank estimation techniques~\cite{Meier21} in the algorithm to determine the appropriate rank and truncating accordingly, without significant additional computation. 

After reviewing preliminary results in Section~\ref{prelim:sec}, we provide an outline of existing randomized algorithms for computing Tucker decomposition in Section~\ref{existingAlgs:sec}. Section~\ref{sec:RTSMS} then describes our algorithm \RTSMS, and we illustrate its performance with experiments in Section~\ref{exper:sec}. 

\section{Preliminaries} \label{prelim:sec}
Let us begin with a brief overview of basic concepts in deterministic Tucker decompositions.

\textbf{Notation} Throughout the paper we denote
$n_{(-k)}:=\Pi_{\substack{j=1\\ j \neq k}}^{d}n_j = \Pi_{j=1}^{d} n_j/n_k$. 

\subsection{Modal unfoldings}
Let $\mathcal{A} \in \mathbb{R}^{n_1 \times n_2 \times \dots \times n_d}$ be a tensor of order $d$. The \emph{mode-$k$ unfolding}
of $\mathcal{A}$, denoted by $A_{(k)}$, is a matrix of size $n_k \times n_{(-k)}$ whose columns are the mode-$k$ fibers of $\mathcal{A}$ \cite[p. 723]{Golubbookori}. 

\subsection{Modal product} A simple but important family of contractions are the modal products. These contractions involve a tensor, a matrix, and a mode. In particular we have the following.
\begin{definition}\label{modalProd:def} 
\cite[p. 727]{Golubbookori}
If $\mathcal{A} \in \mathbb{R}^{n_1 \times n_2 \times \dots \times n_d}$, $M \in \mathbb{R}^{m_k \times n_k}$, and $1 \leq k \leq d$, then \[
\mathcal{B} = \mathcal{A} \times_k M \in\mathbb{R}^{n_1\times \cdots n_{k-1}\times m_k \times n_{k+1}\times\cdots \times n_d}
\] 
denotes the mode-$k$ product of $\mathcal{A}$ and $M$ if $B_{(k)} = M  A_{(k)}$.
\end{definition} 

Note that every mode-$k$ fiber in $\mathcal{C}$ is multiplied by the matrix $M$ requiring $\mbox{size}(\mathcal{C}, k) = \mbox{size}(M, 2)$. Here, $B_{(k)}$ is a matrix of size $m_k \times n_{(-k)}$, hence $\mathcal{B}$ is a tensor of size $n_1 \times \dots n_{k-1} \times m_k \times n_{k+1} \dots \times n_d$.

If $\mathcal{A} \in \mathbb{R}^{n_1 \times n_2 \times \dots \times n_d}$, $F \in \mathbb{R}^{m_1\times n_k}$ and $G \in \mathbb{R}^{m_2\times m_1}$, then 
\begin{equation}
\label{sameModeProdRule:eq}
(\mathcal{A} \times_k F) \times_k G = \mathcal{A} \times_k (G F)
\end{equation}
resulting in a tensor of size $n_1 \times \dots \times n_{k-1} \times m_2 \times n_{k+1} \times \dots \times n_d$. See \cite[p. 461]{kolda2009tensor} or \cite[property 3]{de2000multilinear}.

Reformulating matrix multiplications in terms of modal products, the matrix SVD can be rewritten as follows \cite{Bader06}: 
\begin{equation}
\label{svdMatModal1:eq}
A = U \Sigma V^T \mbox{ means that } A = \Sigma \times_1 U \times_2 V.
\end{equation}

\subsection{Deterministic HOSVD}
For an order-d tensor $\mathcal{A} \in \mathbb{R}^{n_1 \dots \times n_d}$ a HOSVD \cite{de2000multilinear} is a decomposition of the form
\begin{equation}
\label{hosvd1:eq}
\mathcal{A} = \mathcal{C} \times_1 U_1 \times_2 U_2 \times \dots \times_d U_d,
\end{equation}
where the {\em factor matrices} 
$U_k \in \mathbb{R}^{n_k \times r_k}$ $(r_k \leq n_k)$ have orthonormal columns
and the {\em core tensor} 
$\mathcal{C} = \mathcal{A} \times_1 U_1^{T} \times_2 U_2^{T} \dots \times_d U_d^{T} \in \mathbb{R}^{r_1 \times r_2 \dots \times r_d}$
 is not diagonal in general, but enjoys the so-called all-orthogonality property. 
We will use the following standard shorthand notation for Tucker decompositions: 
\[
\mathcal{A} =  \llbracket \mathcal{C}; U_1 , U_2 ,\dots U_d \rrbracket
\]
which, in this case, is defined by the HOSVD \eqnref{hosvd1:eq}.

The Tucker (multilinear) rank of $\mathcal{A}$ is the vector of ranks of modal unfoldings \cite[p. 734]{Golubbookori}
\[
\mbox{rank}(\mathcal{A}) := [r_1, r_2, \dots, r_d]
\]
i.e., $r_k := \mbox{rank}(A_{(k)})$ is the dimension of the space spanned by columns of $A_{(k)}$, which is equal to the span of $U_k$; indeed $U_k$ is computed via the SVD of $A_{(k)}$. 

Note also that the multilinear rank-$(r_1, r_2, \dots, r_d)$ HOSVD truncation of $\mathcal{A}$, while is {\em not} the best multilinear rank-$(r_1, r_2, \dots, r_d)$ approximation to $\mathcal{A}$ in the Frobenius norm, is quasi-optimal; see~\eqnref{HOSVDErrorBnd:eq} below.

\subsection{Deterministic STHOSVD}
This is the sequentially truncated variant of HOSVD which, while processing each mode (in a user-defined processing order $p = [p_1, p_2, \dots, p_d]$, where $1\leq p_k\leq d$ are distinct integers), truncates the tensor. We display its pseudocode in Alg.~\ref{STHOSVD:alg} in supplementary materials.

Let $\hat{\mathcal{C}}^{(i)}$ denote the partially truncated core tensor of size $r_1 \times \dots \times r_i \times n_{i+1} \times \dots \times n_d$ (obtained in the $i$-th step of Alg.~\ref{STHOSVD:alg} where $\hat{\mathcal{C}}^{(0)} := \mathcal{A}$ and $\hat{\mathcal{C}}^{(d)} := \mathcal{C}$). 
Following \cite{Vannieuwenhoven12} we denote the $i$-th partial approximation of $\mathcal{A}$ 
by 
\[
\hat{\mathcal{A}}^{(i)} := \llbracket \hat{\mathcal{C}}^{(i)}; U_1 , \dots, U_i, I, \dots, I \rrbracket = \hat{\mathcal{C}}^{(i)} \times_1 U_1 \times_2 U_2 \dots \times_i U_i,
\]
which has rank-$(r_1, \dots, r_i, n_{i+1}, \dots, n_d)$. In particular we have $\hat{\mathcal{A}}^{(0)} := \mathcal{A}$, and the final approximation obtained by STHOSVD is $\hat{\mathcal{A}}^{(d)} := \hat{\mathcal{A}}$. The factor matrix $U_k$ is computed via the SVD of $\hat A^{(k-1)}_{(k)}$, the mode-$k$ unfolding of $\hat{\mathcal{A}}^{(k-1)}$. This is an overarching theme in the paper; to find an approximate Tucker decomposition, we find a low-rank approximation of the unfolding of the tensor $\hat{\mathcal{A}}^{(k-1)}$.

The following lemma states that the square of the error in approximating $\mathcal{A}$ by STHOSVD is equal to the sum of the square of the errors committed in successive approximations and is bounded by the sum of squares of all the modal singular values that have been discarded\footnote{The statement of the lemma uses a specific processing order but the upper bound remains the same for any other order.}. The error in both STHOSVD and HOSVD satisfy the same upper bound, and STHOSVD tends to give a slightly smaller error~\cite{Vannieuwenhoven12}.

\begin{lemma} \cite{Vannieuwenhoven12, Minster20}
Let $\hat{\mathcal{A}} := \llbracket \mathcal{C}; U_1 , U_2 ,\dots, U_d \rrbracket$ be the rank-$(r_1, r_2, \dots, r_d)$ STHOSVD approximation of $\mathcal{A}$ 
with processing order $p = [1, 2, \dots, d]$ 
Then,
\begin{align*}\label{STHOSVD_error:eq}
\| \mathcal{A} - \hat{\mathcal{A}}\|_F^2 = \sum_{i=1}^{d} \| \hat{\mathcal{A}}^{(i-1)} - \hat{\mathcal{A}}^{(i)} \|_F^2 &\leq 
\sum_{i=1}^{d} \| \mathcal{A} \times_i (I - U_i U_i^T) \|_F^2= \sum_{i=1}^{d} \sum_{j=r_i+1}^{n_i} \sigma_{j}^2(A_{(i)}). \end{align*}
\end{lemma}

Note also that approximations obtained by deterministic HOSVD and STHOSVD are both quasi-optimal in the sense that they satisfy the following bound
\begin{equation}
\label{HOSVDErrorBnd:eq}
\| \mathcal{A} - \hat{\mathcal{A}}\|_F \leq \sqrt{d} \| \mathcal{A} - \hat{\mathcal{A}}_{opt}\|_F,
\end{equation}
where $\hat{\mathcal{A}}_{opt}$ denotes the best rank-$(r_1,r_2, \dots, r_d)$ approximation which can be computed using the HOOI, a computationally expensive nonlinear iteration \cite[pp. 734-735]{Golubbookori}. In many cases, (ST)HOSVD suffices with the near-optimal accuracy~\eqref{HOSVDErrorBnd:eq}.

Our method, RTSMS, is similar to STHOSVD in that in each step it truncates the tensor, but differs from previous methods in how each factor is computed.

\section{Existing randomized algorithms for Tucker decomposition}\label{existingAlgs:sec}
A number of randomized algorithms have been proposed for computing an approximate Tucker decomposition. In what follows we explain some of those techniques with focus on the ones for which implementations are publicly available 
and divide them into two categories: fixed-rank algorithms that require the multilinear rank to be given as an input, and those that  are adaptive in rank; which are the main focus of this paper. 

\subsection{Fixed-rank algorithms} 
Here we suppose that the output target Tucker rank $(r_1,r_2, \dots, r_d)$  is given. 
The first candidate to consider is a randomized analogue of the standard HOSVD algorithm.

\subsubsection{Randomized HOSVD/STHOSVD}

A natural idea to speed up the computation of HOSVD or STHOSVD is to use randomized SVD~\cite{HMT} (Alg.~\ref{RandSVD:alg}) in the computation of the SVD of the unfoldings. This has been done in \cite[Alg. 2]{Zhou14} and \cite[Alg. 3.1, 3.2]{Minster20}, and Minster, Saibaba and Kilmer \cite{Minster20} in particular present extensive analysis of its approximation quality. We display the randomized STHOSVD in Alg.~\ref{Fixed_RSTHOSVD:alg} and compare it against our algorithm in numerical experiments, as it is usually more efficient than the randomized HOSVD. 

\begin{algorithm} 
\caption{Randomized matrix SVD without power iteration 
(Halko-Martinsson-Tropp~\cite{HMT}. See also \cite[Alg. 2.1]{Minster20})
\\ Inputs are matrix $X \in \mathbb{R}^{m \times n}$ and target rank $r$.\\ Output is SVD $X \approx \hat U \hat \Sigma \hat V^T$ with $\hat U \in \mathbb{R}^{m \times r}$, $\hat \Sigma \in \mathbb{R}^{r \times r}$, and $\hat V \in \mathbb{R}^{n \times r}$. \label{RandSVD:alg}}
\begin{algorithmic}[1]
\STATE Draw a sketch (e.g. Gaussian) matrix $\Omega \in \mathbb{R}^{n \times \hat r}$, and compute $Y := X \Omega$. \hfill  
\STATE Compute thin QR decomposition $Y := Q R$. \hfill  
     \COMMENT{$Q$: approximate range$(X)$}
\STATE Compute $B := Q^T\! X$. \STATE Compute thin SVD $B := \hat U_B \hat \Sigma \ \! \hat V^T$.\hfill  \COMMENT{$\hat U_B$ and $\hat \Sigma$ are $\hat r \times \hat r$ while $\hat V$ is $n \times \hat r$.} \STATE Output $\hat U := Q\ \!  \hat U_B(:, 1:r)$, $\hat \Sigma := \hat \Sigma(1:r, 1:r)$, and 
$\hat V := \hat V(:, 1:r)$.
\end{algorithmic}
\end{algorithm}

A well-known technique to improve the quality of the approximation obtained by randomized matrix SVD is to run a few iterations of power method~\cite{HMT}. Mathematically, this means $Y$ in Alg.~\ref{RandSVD:alg} is replaced with $X(X^TX)^q$ for an integer $q$. \begin{algorithm}
\caption{Randomized STHOSVD \cite[Alg. 3.2]{Minster20}. \\ Inputs are $\mathcal{A} \in \mathbb{R}^{n_1 \times n_2 \times \dots \times n_d}$, multilinear rank $(r_1, r_2, \dots, r_d)$, oversampling parameter $\tilde p$, and processing order $\bf{p}$ of the modes.\\ Output is $\mathcal{A} \approx   \llbracket \mathcal{C}; U_1 , U_2 ,\dots, U_d \rrbracket$.}\label{Fixed_RSTHOSVD:alg}
\begin{algorithmic}[1]
\STATE Set $\mathcal{C}:= \mathcal{A}$.
\FOR{$i=p_1,\ldots,p_d$}
    \STATE Draw a sketch (e.g. Gaussian) matrix  $\Omega_i$ of size $z_{i} \times \hat r_{i}$ where $\hat r_{i}:= r_{i} + \tilde p$.
\STATE Set $U_i$ as left singular vectors computed with $\texttt{RandSVD}(C_{(i)}, r_i, \Omega_i)$ using Alg.~\ref{RandSVD:alg}. 
\STATE Update $C_{(i)} := U_i^T C_{(i)}$. \hfill       \COMMENT{Overwriting $C_{(i)}$ overwrites $\mathcal{C}$.}
\ENDFOR
\end{algorithmic}
\end{algorithm}

\paragraph{Note on sketching}
An important aspect of any randomized algorithm utilizing a random sketch is: which sketch should be used? In \cite{Minster20} and \cite[Alg. 2]{Zhou14}, they are taken to be Gaussian matrices. This class of random matrices usually come with the strongest theoretical guarantees and robustness. Other more structured random sketches have been proposed and shown to be effective, including sparse~\cite{woodruff2014sketching} and FFT-based~\cite{Tropp2011ImprovedTransform} sketches, and in the tensor sketching context, those employing Khatri-Rao products~\cite{Sun20}. 

In this paper we mostly focus on Gaussian sketches, unless otherwise specified---while other sketches may well be more efficient in theory, the smallness of our sketches (an important feature of our algorithm RTSMS) means that the advantages offered by structured sketches are likely to be limited. 

\paragraph{The one-pass algorithms of Malik and Becker}\label{TS:subsec}
Tucker-TS and Tucker-TTMTS are two algorithms for Tucker decomposition \cite[Alg. 2]{Malik18} which rely on the 
TensorSketch framework \cite{Pagh13} and can handle tensors whose elements are streamed, i.e., they are lost once processed. These algorithms are variants of the standard ALS (HOOI) and iteratively employ sketching for computing solutions to large overdetermined least-squares problems, and for efficiently approximating chains of tensor-matrix products which are Kronecker products of smaller matrices. Inherited from the characteristics of the TensorSketch, these algorithms only require a single pass of the input tensor.

\paragraph{The one-pass algorithm of Sun, Guo, Luo, Tropp and Udell}\label{singlePass:subsec}
Another single-pass sketching algorithm targeting streaming Tucker decomposition is introduced in \cite{Sun20} where a rigorous theoretical guarantee on the approximation error was also provided. Treating the tensor as a multilinear operator, it employs the Khatri-Rao product of random matrices to identify a low-dimensional subspace for each mode of the tensor that captures the action of the operator along that mode, and then produces a low-rank operator with the same action on the identified low-dimensional tensor product space, which is helpful especially when storage cost is a potential bottleneck.

\subsection{Rank-adaptive algorithms}
We now turn to existing algorithms that are rank-adaptive, i.e., they are able to determine the appropriate rank on the fly to achieve a prescribed approximation tolerance. Although the first natural candidate to consider is the randomized adaptive variant of the standard HOSVD algorithm \cite[Alg. 4.1]{Minster20}, in what follows we focus on its sequentially truncated version \cite[Alg. 4.2]{Minster20} which is more efficient in practice. 

\subsubsection{Adaptive R-STHOSVD}
We outline the adaptive randomized STHOSVD method \cite[Alg. 4.2]{Minster20} below in Alg.~\ref{adapt_RST:alg} for the sake of completeness. 

\begin{algorithm} 
\caption{Adaptive R-STHOSVD \cite[Alg. 4.2]{Minster20}. \\ Inputs are $\mathcal{A} \in \mathbb{R}^{n_1 \times n_2 \times \dots \times n_d}$, tolerance {\em tol} $\geq 0$, blocking parameter $b \geq 1$ and processing order $\bf{p}$ of the modes.\\ Output is $\mathcal{A} \approx   \llbracket \mathcal{C}; \hat U_1 , \hat U_2 ,\dots, \hat U_d \rrbracket$.}
\label{adapt_RST:alg}
\begin{algorithmic}[1]
\STATE Set $\mathcal{C}:= \mathcal{A}$.
\FOR{$i=p_1,\ldots,p_d$}
\STATE Compute $\hat U_i = \texttt{AdaptRandRangeFinder}(C_{(i)}, \frac{{\rm tol}}{\sqrt{d}}, b)$.\hfill \COMMENT{Error-controlled range finder. We use \texttt{svdsketch}.}
    \STATE Update $C_{(i)} = \hat U_i^T C_{(i)}$. \hfill       \COMMENT{Overwriting $C_{(i)}$ overwrites $\mathcal{C}$.}
\ENDFOR
\STATE Return $\mathcal{C}$ by tensorizing $C_{(i)}$.
\end{algorithmic}
\end{algorithm}

At the heart of Alg.~\ref{adapt_RST:alg} is an adaptive randomized matrix range finder. Numerous techniques have been introduced in the literature for this task \cite{HMT, Martinsson16, Yu18}. Given a matrix $A$ and a positive tolerance tol, the primary objective is to find a (tall) matrix $Q$ with orthonormal columns such that the relative residual in approximating the range of $A$ by $Q$ is bounded from above by tol, i.e., 
\begin{equation}
\label{adaptRangFinderCond:eq}
\| A - Q Q^T A \| \leq {\rm tol} \| A\|.
\end{equation}
The rank of the low-rank approximation then corresponds to the number of columns in $Q$. The idea of adaptive randomized range finders is to start with a limited number of columns in the random matrix $\Omega$ in order to estimate the range $Q$. Then, the number of random columns drawn is gradually increased until $Q$ satisfies \eqnref{adaptRangFinderCond:eq}. Among the current state-of-the-art randomized rangefinders is the randQB\_EI\_auto algorithm of \cite{Yu18}, which is a variant of an algorithm originally introduced in \cite{HMT}. randQB\_EI\_auto has been integrated into MATLAB  since release 2020b, as the \texttt{svdsketch} function. 
We will compare RTSMS with Alg.~\ref{adapt_RST:alg} whose Step 4 calls \texttt{svdsketch}. See Section~\ref{exper:sec} for more details.

\section{\RTSMS}\label{sec:RTSMS}
We now describe \RTSMS, our new randomized algorithm for computing a Tucker decomposition and HOSVD. \RTSMS\ follows the basic structure of R-STHOSVD of sequentially finding a low-rank approximation to the unfoldings, but crucially, it only applies a random sketch on \emph{one mode} at each step, resulting in a dramatically smaller sketch matrix. \RTSMS\ can thus offer speedup over existing algorithms, especially when sketching the large dimensions is costly.

\RTSMS\ is shown in pseudocode in Alg.~\ref{RTSMS:alg}. In a nutshell, 
we find low-rank approximations to the unfolding matrices $A_{(i)}$ (or more precisely the unfolding of the current core tensor $\mathcal{B}^{\rm{old}}$) via a single-mode sketch. 
Namely, to approximate $A_{(1)}\in\mathbb{R}^{n_1\times n_{(-1)}}$, we sketch from the left by an operation equivalent to computing $\Omega_1A_{(1)}$ where $\Omega_1 \in \mathbb{R}^{\hat r_1 \times n_1}$, and then find  $F_1\in\mathbb{R}^{n_1\times \hat r_1}$
such that 
$F_1\Omega_1A_{(1)}\approx A_{(1)}$
 by solving a massively overdetermined least-squares problem with many right-hand sides $\min_{F_1}\|
A_{(1)}^T\Omega_1^TF_1^T-  A_{(1)}^T\|_F$, which we do by randomized row subset selection, regularization, and iterative refinement. We discuss the details of the least-squares solution in Section~\ref{alg:LS}.

Another key component of \RTSMS\ is the rank adaptivity. In Steps 6-19, we find an appropriate Tucker rank on the fly by employing a fast rank estimator~\cite[Alg.~1]{Meier21} on each unfolding, which are truncated in turn. 
The rank is estimated more efficiently than \texttt{svdsketch}, and much of the computation needed for rank estimation can be reused within \RTSMS; for example, the first $\min(\tilde r_i,\hat r_i)$ rows of $\Omega_i,\Omega_iM$ in steps 21 and 22 are identical, and the QR factorization $(\Omega_i M Y_i)^T=QR$ can be reused by Alg.~\ref{RTSMS:LSalg}.
Alternatively, \RTSMS\ is also able to take the rank as a user-defined input.

To illustrate the high-level ideas, let us provide an overview of \RTSMS\ when applied to compute a rank-$(r_1, r_2, r_3)$ HOSVD of a tensor $\mathcal{A}$ of size $n_1 \times n_2 \times n_3$ with the processing order of ${\bf p} = [1\ 2\ 3]$. We illustrate the process in Figure~\ref{tucker:fig}. In Alg.~\ref{RTSMS:alg} a Tucker decomposition $\mathcal{A} \approx \llbracket \mathcal{B}^{\rm{new}}; F_1 , F_2, F_3 \rrbracket$ is computed where the $i$-th factor matrix $F_i$ is of size $n_i \times \hat r_i$ (where $\hat r_i\approx 1.5r_i$), rather than $n_i \times r_i$, does not have orthonormal columns in general, and the size of the final core tensor $\mathcal{B}^{\rm{new}}$ is $(\hat r_1, \hat r_2, \hat r_3)$. More specifically,
\begin{enumerate}
\item At the beginning, $i = p_1 = 1$ and the overall goal is to compute a Tucker1 decomposition \cite[\S~4.5.1]{Kroonenberg08}
\[
\mathcal{A} \approx \mathcal{B}^{\rm{new}} \times_1 F_1.
\]
See Figure~\ref{tucker:fig}, top. Importantly, we directly compute the ``temporary core tensor'' $\mathcal{B}^{\rm{new}} =\mathcal{A}\times_1\Omega_1$ as the sketch of $\mathcal{A}$ in Steps 20-21, then find the factor matrix 
$F_1$ in Step 22 of Alg.~\ref{RTSMS:alg}. This is in contrast to existing approaches, e.g. (R-)STHOSVD, where $F_1$ is computed and orthonormalized first and then $\mathcal{B}^{\rm{new}}$  is taken to be $\mathcal{A}\times_1 F_1^T$; see Step 5 of Alg.~\ref{Fixed_RSTHOSVD:alg}. This means that \RTSMS\ is almost single-pass, i.e., it does not need to revisit $\mathcal{A}$ once the sketch $\mathcal{A}\times_1\Omega_1$ is computed, aside from the need to subsample a small number of fibers of $\mathcal{A}$, 
see Section~\ref{sec:LSsolvesec}. 

\begin{figure}[!h]
\begin{center}
\includegraphics[width=0.8\textwidth]{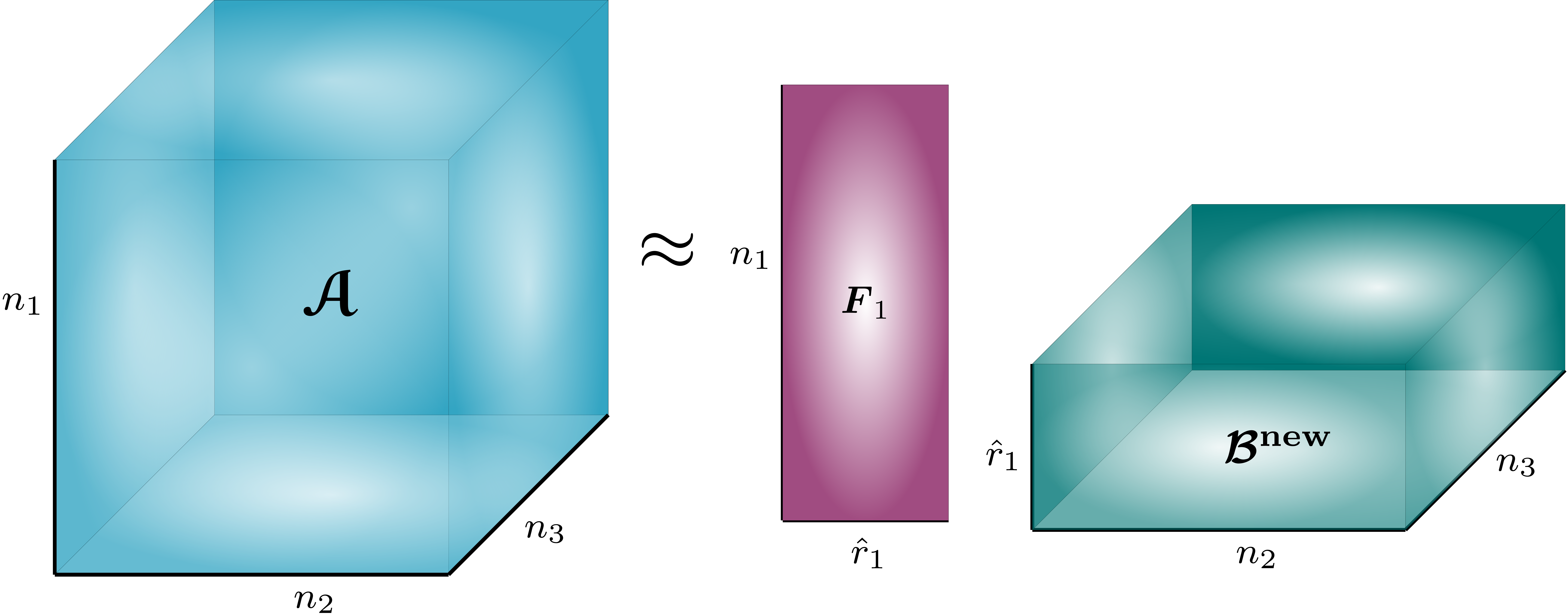}
\includegraphics[width=0.8\textwidth]{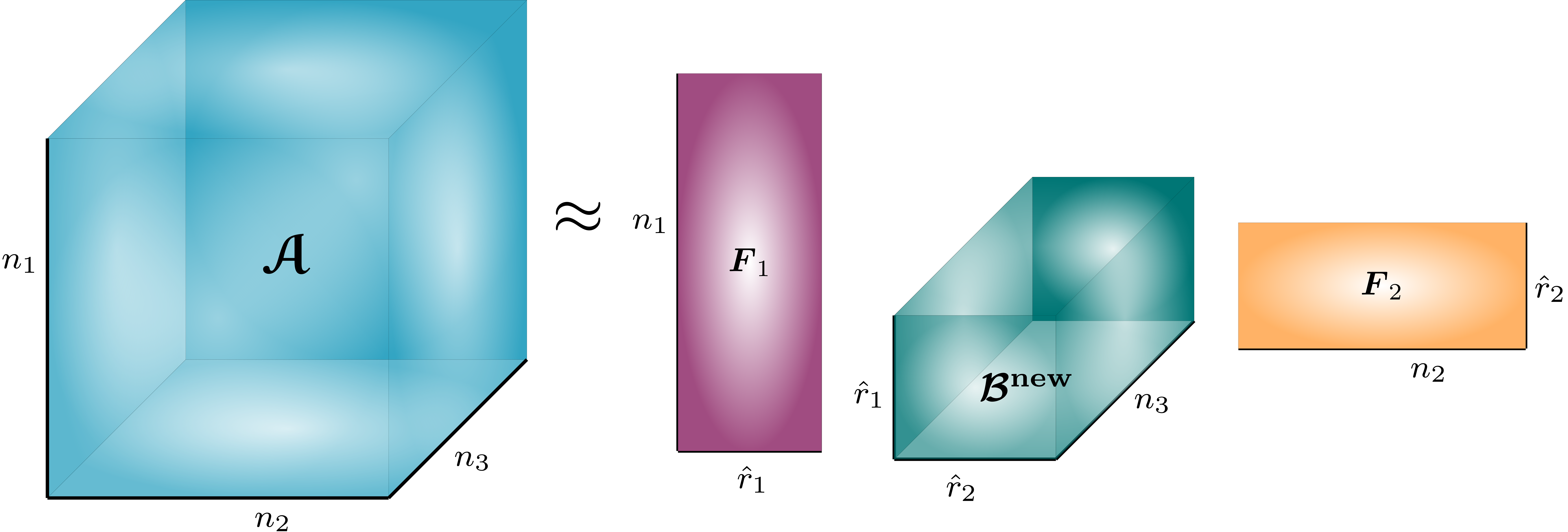}
\includegraphics[width=0.8\textwidth]{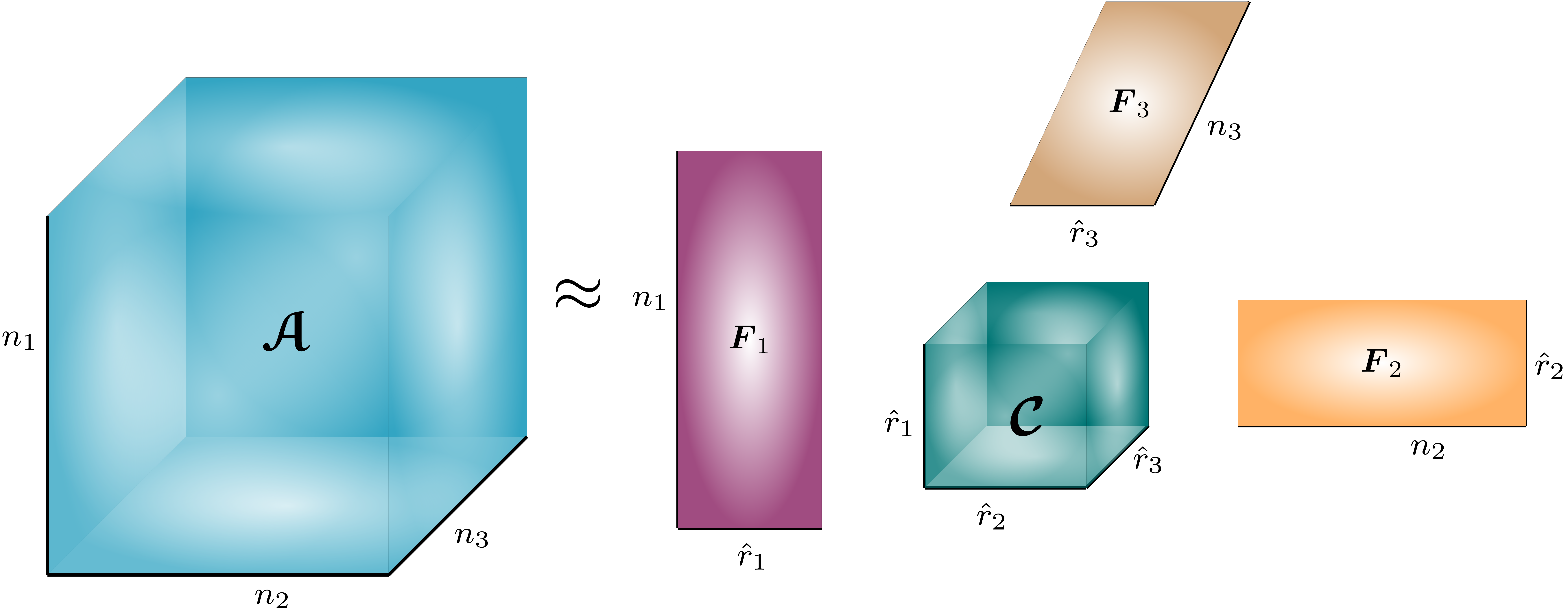}
\caption{Illustration of RTSMS for 3D tensors, Tucker rank $(r_1, r_2, r_3)$ and mode-processing order $[1\ 2\ 3]$. Top: Tucker1 decomposition $\mathcal{A} \approx \mathcal{B}^{\rm{new}} \times_1 F_1$ where $\mathcal{A} =: \mathcal{B}^{\rm{old}} \in \mathbb{R}^{n_1 \times n_2 \times n_3}$, $\mathcal{B}^{\rm{new}} \in \mathbb{R}^{\hat r_1 \times n_2 \times n_3}$ and $F_1 \in \mathbb{R}^{n_1 \times \hat r_1}$. In \RTSMS, $\mathcal{B}^{\rm{new}}$ is computed before $F_1$.
Middle: Tucker2 decomposition $\mathcal{A} \approx \mathcal{B}^{\rm{new}} \times_1 F_1 \times_2 F_2$ where $\mathcal{B}^{\rm{new}} \in \mathbb{R}^{\hat r_1 \times \hat r_2 \times n_3}$ and $F_2 \in \mathbb{R}^{n_2 \times \hat r_2}$. Note that $\mathcal{B}^{\rm{new}}$ here is overwritten on $\mathcal{B}^{\rm{new}}$ from top and therefore represents a different tensor. Bottom: Tucker3 decomposition 
$\mathcal{A} \approx \mathcal{B}^{\rm{new}} \times_1 F_1 \times_2 F_2 \times_3 F_3$ where this time $\mathcal{B}^{\rm{new}} \in \mathbb{R}^{\hat r_1 \times \hat r_2 \times \hat r_3}$ is the final core tensor $\mathcal{C}$ computed once Step 24 of Alg.~\ref{RTSMS:alg} is executed and $F_3 \in \mathbb{R}^{n_3 \times \hat r_3}$. Figure generated using {\tt Tensor$\_$tikz}~\cite{Kour23}.}
\label{tucker:fig}
\end{center}
\end{figure}

\begin{itemize}
\item 
In Steps 20-21 we 
find an approximate row space of $A_{(1)}$ using a randomized sketch. We first compute the mode-1 product of $\mathcal{A}$ and $\Omega_1$, which is equivalent to $B_{(1)}^{\rm{new}} = \Omega_1 B_{(1)}^{\rm{old}}$. Here, $B_{(1)}^{\rm{old}} = A_{(1)}$ is a matrix of size $ n_1 \times (n_2 n_3)$, $\Omega_1$ is $\hat r_1 \times n_1$ and so $\mathcal{B}^{\rm{new}}$ is a tensor of size $\hat r_1 \times n_2 \times n_3$, where $\hat r_1$ is found by the rank estimator. The cost of computing $\mathcal{B}^{\rm{new}}$ is $\mathcal{O}(n_1 n_2 n_3 \hat r_1)$ which is cubic in the large dimensions. 
Since $\Omega_1$ is random by construction, provided that $\sigma_{r_1}(A_{(1)})<tol$, with high-probability 
the row space of $B_{(1)}^{\rm{new}}$ is a good approximation to the dominant row space of $A_{(1)}$.

\item The theory in randomized low-rank approximation indicates that there exists 
a factor matrix $F_1$ such that 
$\mathcal{B}^{\rm{new}} \times_1 F_1\approx \mathcal{B}^{\rm{old}}$. In Step 22 we aim at solving the least-squares problem 
  \begin{equation}
    \label{eq:lsfirst}
\min_{F_1}\|\mathcal{B}^{\rm{new}} \times_1 F_1 - \mathcal{B}^{\rm{old}}\|_F.    
  \end{equation}
Its efficient and robust solution proves to be subtle, and we defer the discussion to Section~\ref{alg:LS}.
\end{itemize}

\item Next $i = p_2 = 2$ and we aim at decomposing the updated $\mathcal{B}^{\rm{old}}$ from the previous step hence eventually computing a Tucker2 decomposition of the original tensor $\mathcal{A}$. See \cite[chap. 4.5.2]{Kroonenberg08} and the middle picture of Figure~\ref{tucker:fig}. This is done by essentially repeating the process in the previous step on the second mode of the temporary core tensor $\mathcal{B}^{\rm{old}}$:

\begin{itemize}
\item In Steps 20-21 we 
sketch to find the row space of $B_{(2)}^{\rm{old}}$. We first compute the mode-2 product of $\mathcal{B}^{\rm{old}}$ and $\Omega_2$ which is equivalent to $B_{(2)}^{\rm{new}} = \Omega_2 B_{(2)}^{\rm{old}}$. Here, $B_{(2)}^{\rm{old}}$ is a matrix of size $n_2 \times (n_3 \hat r_1)$, $\Omega_2$ is $\hat r_2 \times n_2$ and so the updated temporary core tensor $\mathcal{B}^{\rm{new}}$ is a tensor of size $\hat r_1 \times \hat r_2 \times n_3$. This process shrinks the dimension of the row space of $A_{(2)}$ from the large $n_2$ to the smaller $\hat r_2$ and is the first task 
in computing a decomposition of the temporary core tensor, and equivalently, a Tucker2 decomposition of $\mathcal{A}$. As before, the row space of $B_{(2)}^{\rm{new}}$ approximates that of $A_{(2)}$ with high probability, provided $\hat r_2$ was found appropriately by the rank estimator.\item Step 22 then solves the least-squares problem
$
\min_{F_2}\|\mathcal{B}^{\rm{new}} \times_2 F_2 - \mathcal{B}^{\rm{old}}\|_F.    $

\end{itemize}

\item Finally, $i = p_3 = 3$ and we aim at decomposing the updated $\mathcal{B}^{\rm{old}}$ from the previous step hence eventually computing a Tucker3 (or simply a Tucker) decomposition of the original tensor $\mathcal{A}$. See \cite[chap. 4.5.3]{Kroonenberg08} and Figure~\ref{tucker:fig}, bottom (adapted from \cite{Kour23}).

\end{enumerate}

The computational complexity of the entire algorithm is dominated by the first tensor-matrix multiplication $\mathcal{A} \times_1 \Omega_1$, which is equivalent to 
computing $\Omega_1A_{(1)}$. 

\paragraph{Single-mode vs. multi-mode sketch}
 As far as we know, existing randomized algorithms for Tucker decomposition require applying dimension reduction maps to the tensor unfoldings from the {\em right} (i.e., the larger dimension). Focusing on mode-1, one faces $A_{(1)}$ which is a short fat matrix of size $n_{1} \times n_{(-1)}$. Sketching from the right then involves computing with a tall skinny randomized matrix, whose size is $n_{(-1)} \times \hat r_1$. R-STHOSVD and R-HOSVD are among the examples, as are most algorithms we mentioned. Considering the fact that $n_{(-1)}$ could be huge, generating or even storing the randomized sketch matrix could be problematic. The focus of the single-pass techniques \cite{Malik18,Sun20}
mentioned in section~\ref{TS:subsec} has therefore been on designing algorithms which do not need generating and storing a tall skinny dimension reduction map. Instead, the fat side of the unfolding is sketched by skillfully employing Khatri-Rao products of smaller random matrices so as to avoid explicitly forming random matrices which have $n_{(-1)}$ rows. In this sense, one can classify all the aforementioned algorithms as those which sketch in all-but-one modes when processing each of the modes.

In Alg.~\ref{RTSMS:alg}, the original tensor $\mathcal{A}$, which is potentially large in all dimensions, is directly used only in the first recursion $i=p_1$. The updating of $\mathcal{B}^{\rm{old}}$ in Step 23 is where sequential truncation takes place, making sure that subsequent computations proceed with the potentially much smaller truncated tensors. 
This truncation can improve the efficiency significantly, just as STHOSVD does over HOSVD.
Notice however that the output Tucker rank will be ${\bf \hat r}$, which is slightly larger than ${\bf r}$.

\begin{algorithm}[!h] 
\caption{\texttt{RTSMS}: Randomized Tucker with single-mode sketching \\ Inputs are $\mathcal{A} \in \mathbb{R}^{n_1 \times n_2 \times \dots \times n_d}$, target tolerance {\em tol} on the relative residual, and processing order $\bf{p}$ of the modes.\\ Output is Tucker decomposition $\mathcal{A} \approx \llbracket \mathcal{C}; F_1 , F_2 ,\dots, F_d \rrbracket$.}
\label{RTSMS:alg} \begin{algorithmic}[1]
\STATE Initialize ${\bf r}$ with a small ambitious rank estimate like $(10, 10, \dots, 10)$.
\STATE Set $\mathcal{B}^{\rm{old}} := \mathcal{A}$.
\FOR{$i=p_1,\ldots,p_d$}
\STATE Set ${\rm flag}_{i}$ to `unhappy'. 
\STATE Find $M$, the $i$-th modal unfolding of $\mathcal{B}^{\rm{old}}$.
\WHILE{${\rm flag}_{i}$ is `unhappy'} \hfill \COMMENT{rank estimation}
    \STATE Set $\tilde r_{i} := {{\rm round}(1.1\ r_{i})}$.
    \STATE Draw a standard random Gaussian matrix $\Omega_i$ of size $\tilde r_{i} \times n_{i}$.
    \STATE Form the matrix $\Omega_i M$ of size $\tilde r_{i} \times z_{i}$ where $z_i := (\Pi_{j=1}^{i-1} \tilde r_j) (\Pi_{j=i+1}^{d} n_j)$.
    \STATE Set $s_{i} := k  \tilde r_i$ (e.g. $k=4$), draw SRFT 
$Y_i\in\mathbb{R}^{z_i \times s_i}$ 
and form $\Omega_i M Y_i\in\mathbb{R}^{\tilde r_i \times s_i}$.
\STATE Compute thin QR factorization $(\Omega_i M Y_i)^T=QR$. 
\STATE Find smallest $\ell$ such that
     $\sigma_{\ell+1} (R) \leq {\rm tol}\ \sigma_{1} (R)$.
\IF{$\ell < r_{i}$} 
\STATE Change ${\rm flag}_{i}$ to `happy', set $i$-th multilinear rank estimate $r_i$ to $\ell$.
\STATE Update $\mathcal{B}^{\rm{old}}$ to be $\Omega_i M$ tensorized in the $i$-th mode.
\ELSIF{$\ell$ is `empty' or $\ell = \tilde r_{i}$ 
implying the sketch was full-rank}
\STATE Increase the rank estimate, e.g. $r_{i}:={{\rm round}(1.7 r_{i})}$, return to Step 7.
\ENDIF
\ENDWHILE 
\STATE Let $\Omega_i$ be $\hat r_{i} \times n_{i} $ Gaussian, where $\hat r_{i}:= {{\rm round}(1.5\ r_{i})}$.
\hfill       \COMMENT{low-rank approx}
    \STATE Compute $\mathcal{B}^{\rm{new}} = \mathcal{B}^{\rm{old}} \times_i \Omega_i$ (reusing $\Omega_iM$ from step 9)
\STATE Find $F_i$ of size $n_i \times \hat r_i$ to minimize $\| \mathcal{B}^{\rm{new}} \times_i F_i - \mathcal{B}^{\rm{old}} \|_F$, using Alg.~\ref{RTSMS:LSalg} 
     \STATE Update $\mathcal{B}^{\rm{old}} := \mathcal{B}^{\rm{new}}$.
\ENDFOR
\STATE Set $\mathcal{C} := \mathcal{B}^{\rm{new}}$.
\end{algorithmic}
\end{algorithm}

\paragraph{Parallel, non-sequential variant} We have emphasized the advantages of \RTSMS, most notably the fact that the sketches are small. 
One possible advantage of a non-sequential algorithm is that they are highly parallalizable and suitable for distributed computation, as dicussed e.g. in \cite[p. 6]{Zhou14}.
It is possible to design a variant of \RTSMS\ that works parallely: run the first step ($i=p_1$ in Alg.~\ref{RTSMS:alg}) to find the factor matrices $F_i$ for each $i$ from $\mathcal{A}$.
We then find the core tensor as follows: compute the thin QR factorizations $F_i=U_iR_i$ for each $i$, and project them onto $\mathcal{A}$, i.e., $\mathcal{C}=\mathcal{A}\times_1 U_1^T\cdots \times_d U_d^T$. 
We do not discuss this further, as in our sequential experiments the standard \RTSMS\ is more efficient with a lower computational cost.

\subsection{Analysis of \RTSMS}
Let us explain why \RTSMS\ is able to find an approximate Tucker decomposition. We focus on the first step $i=1$ (and assume WLOG $p_1=1$), as the other cases are essentially a repetition. 

Suppose that $\mathcal{A}$ has an approximate HOSVD 
(here we assume the factors 
$F_i$ are orthonormal and denote them by $U_i$, which simplifies the theory but is not necessary in the algorithm)
\begin{equation}  \label{eq:approxhosvd}
\mathcal{A} = \mathcal{C} \times_1 U_1 \times_2 U_2 \times \dots \times_d U_d + \mathcal{E},  
\end{equation}
where 
each $U_k\in\mathbb{R}^{n_k\times r_k}$ 
is orthonormal 
$U_k^TU_k=I_{r_k}$, and $\mathcal{E}$ is small in norm; that is, $\mathcal{A}$ has an approximate Tucker decomposition of rank $(r_1,r_2, \dots,r_d)$. 
Then in the first step of the algorithm
$\mathcal{B}^{\rm{new}} = \mathcal{A} \times_1 \Omega_1$, one obtains 
$\mathcal{B}^{\rm{new}} = \mathcal{C} \times_1 (\Omega_1U_1) \times_2 U_2 \times \dots \times_d U_d + \mathcal{E}\times_1 \Omega_1.$ 
Note that $\|\mathcal{E}\times_1 \Omega_1\|_F\leq O(\|\mathcal{E}\|_F)$, because multiplication by Gaussian matrices roughly preserves the norm~\cite[\S~10]{HMT} (the $O$ notation hides constant multiples of $\sqrt{n_1}$).
Now since $\Omega_1$ is Gaussian so is $\Omega_1 U_1$ by orthogonal invariance, and it is a (tall) $\hat r_1\times r_1$ rectangular Gaussian matrix; therefore well-conditioned with high probability by the Marchenko-Pastur rule~\cite{pastur1967distribution}, or more specifically Davidson-Szarek's result~\cite[Thm.~II.13]{davidson2001local}. 

In terms of the mode-1 unfolding, we have 
$(\mathcal{A}\times_1\Omega_1)_{(1)} = \Omega_1 A_{(1)}$. 
Now recalling~\eqref{eq:approxhosvd} note that the unfolding of 
$\mathcal{C} \times_1 U_1 \times_2 U_2 \times \dots \times_d U_d$ can be written as $U_1G_1$ for some $G_1\in\mathbb{R}^{r_1\times n_{(-1)}}$, so by assumption the mode-1 unfolding of $\mathcal{A}$ is $A_{(1)} = U_1G_1 + \tilde E$, where $\|\tilde E\|_F=\|\mathcal{E}\|_F$. As $U_1$ is $n_1\times r_1$, this implies that the matrix $A_{(1)}$ can be approximated in the Frobenius norm by a rank-$r_1$ matrix up to $\|\tilde E\|_F$.

This is precisely the situation where randomized algorithms for low-rank approximation are highly effective. In particular by the analysis in \cite[\S~10]{HMT} (applied to $A_{(1)}^T$ rather than $A_{(1)}$), it follows that by taking $\Omega_1$ to have $\hat r_1>r_1$ rows (say $\hat r_1=1.2r_1$), the row space of $\Omega_1 A_{(1)}$ captures that of $A_{(1)}$ up to a small multiple of $\|\tilde E\|_F$. This implies that using the thin QR factorization $(\Omega_1 A_{(1)})^T=QR$, the rank-$\hat r_1$ matrix 
\begin{equation}
  \label{eq:HMTapprox}
A_{(1)}QQ^T  \approx A_{(1)}
\end{equation}
 approximates $A_{(1)}$ up to a modest multiple of $\|\tilde E\|_F$, and hence of $\|\mathcal{E}\|_F$.  

Note that this discussion shows that $\min_{U_1\in\mathbb{R}^{n_1\times \hat r_1}}\| U_1 (\Omega_1A_{(1)}) -A_{(1)}\|_F=O(\|\mathcal{E}\|_F)$, because $Q^T$ and $\Omega_1A_{(1)}$ have the same row space. 
This is equivalent to the least-squares problem~\eqref{eq:lsfirst}.

\subsection{Algorithm to $\mathbf{\minimizebf_{\hat F_i}\|  \hat F_i(\Omega_iA_{(i)}) -A_{(i)}\|_F}$}
\label{alg:LS} 
In RTSMS we do not form $Q$ or $A_{(1)}QQ^T $, as these operations can be expensive and dominate the computation. 
In particular, the computation of $A_{(1)}Q$ involves the large dimension $n_{(-1)}$ and can be expensive. 
Instead, in \RTSMS\ we attempt to directly find the factor matrix $\hat F_1\in\mathbb{R}^{n_1\times \hat r_1}$ via minimizing $\| \hat F_1 (\Omega_1A_{(1)}) -A_{(1)}\|_F$, which we rewrite in standard form of a least-squares (LS) problem as 
\begin{equation}  \label{eq:LSstate}
\minimize_{\hat F_1\in\mathbb{R}^{n_1 \times \hat r_{1}}}\|(A_{(1)}^T\Omega_1^T) \hat F_1^T - A_{(1)}^T\|_F.
\end{equation}
This LS problem has several important features worth noting: (i) it is massively overdetermined $A_{(1)}^T\Omega_1^T\in \mathbb{R}^{n_{(-1)}\times \hat r_1}$, (ii) it has many ($n_1$) right-hand sides, and (iii) the coefficient matrix $A_{(1)}^T\Omega_1^T$ is ill-conditioned, and numerically rank-deficient (by design, assuming $tol$ is close to machine precision).\ Solving~\eqref{eq:LSstate} exactly via the classic QR-based approach gives the approximation $QQ^TA_{(1)}^T\approx A_{(1)}^T$, which is equivalent to~\eqref{eq:HMTapprox}. We employ three techniques to devise a more efficient (yet robust) approximation algorithm.

\subsubsection{Randomized sketching}
In order to speed up the computation, we solve the LS problem~\eqref{eq:LSstate} using randomization. This is now a standard technique for solving highly-overdetermined least-squares problem, of which \eqref{eq:LSstate} is one good example. 

Among the most successful ideas in the randomized solution of LS problems is \emph{sketching}, wherein instead of~\eqref{eq:LSstate} one solves the sketched problem
\begin{equation}  \label{eq:LSsketch}
\minimize_{\hat F_1\in\mathbb{R}^{n_1 \times \hat r_{1}}}\|S(A_{(1)}^T\Omega_1^T \hat F_1^T - A_{(1)}^T)\|_F,
\end{equation}
where $S\in\mathbb{R}^{s_1 \times n_{(-1)}}$ ($s_1 \geq  \hat r_1$) is a random matrix, called the sketching matrix. Effective choices of $S$ include Gaussian, FFT-based (e.g. SRFT) sketches~\cite{Tropp2011ImprovedTransform} and sparse sketches~\cite{clarkson2017low}. The solution for~\eqref{eq:LSsketch} is $\hat F_1^T= (SA_{(1)}^T\Omega_1^T)^\dagger (SA_{(1)}^T)$.

The choice of $S$ that is the easiest to analyze is when it is taken to be a random Gaussian matrix. Then 
the resulting rank-$\hat r_1$ approximation $A_{(1)}\approx \hat F_1 \Omega_1 A_{(1)}$ obtained by solving~\eqref{eq:LSsketch} becomes equal to the generalized Nystr\"om (GN) approximation~\cite{woolfe2008fast,clarkson2017low,Nakatsukasa20fast}, which takes  $ A_{(1)}\approx  A_{(1)}X(YA_{(1)}X)^\dagger YA_{(1)}$, where $X,Y$ are random sketches of appropriate sizes. To see this, note that with the solution of~\eqref{eq:LSsketch} $\hat F_1=A_{(1)}S^T(\Omega_1 A_{(1)}S^T)^\dagger$ we have $\hat F_1 \Omega_1 A_{(1)} =A_{(1)}S^T(\Omega_1 A_{(1)}S^T)^\dagger \Omega_1 A_{(1)}$; i.e., they coincide by taking $S^T=X$ and $\Omega_1=Y$. 

Despite the close connection, an important difference here is that in \RTSMS, $S^T=X$ will be generated using $\Omega_1 A_{(1)}$ as we describe in Section~\ref{sec:LSsolvesec}, so it depends on $\Omega_1=Y$, unlike the standard GN where $X,Y$ are independent. Furthermore, taking $X$ to be an independent sketch necessitates sketching $A_{(1)}$ from the right, which violates our `single-mode-sketch' approach and results in inefficiency. 

Here we establish a result on the accuracy of the approximate Tucker decomposition obtained in the GN way. 
While this is not exactly what our algorithm does (which uses a different sketch and solves~\eqref{eq:LSsketch} using techniques including regularization and iterative refinement as we describe shortly), 
we record it as it is known that analyses of Gaussian sketches tend to reflect the typical behavior of algorithms that employ randomized sketching quite well~\cite{MartinssonTroppacta}, and the techniques described in Section~\ref{sec:LSsolvesec} are essentially attempts at solving the ill-conditioned problem~\eqref{eq:LSsketch} in a numerically stable fashion.
\begin{theorem}
Let $\mathcal{\hat A}:=\llbracket \mathcal{C}; \hat F_1 , \hat F_2 ,\dots, \hat F_d \rrbracket$ 
be the output of 
RTSMS (Algorithm~\ref{RTSMS:alg}), where the least-squares problems in step 22 are solved via~\eqref{eq:LSsketch}, where 
each $S$ is taken\footnote{The sketch $S$ of course depends on the step $i$, so perhaps should be denoted by $S_i$; we drop the subscript for simplicity and consistency with the remainder.} to be 
$s_i \times z_i$
Gaussian. Then 
\begin{equation}
  \label{eq:theorembound}
\mathbb{E}\|\mathcal{\hat A}-\mathcal{A}\|_F\leq 
\sum_{j=1}^d
\left(\prod_{i=1}^j 
\sqrt{1+\frac{s_i}{s_i-\hat r_i-1}}\sqrt{1+\frac{\hat r_i}{\hat r_i-\ell_i-1}}\right)\|\mathcal{A}-\mathcal{A}_{{\rm opt}}\|_F,   
\end{equation}
where $\mathcal{A}_{{\rm opt}}$ is the best Tucker approximation  of rank $(r_1, r_2, \dots,r_d)$ to $\mathcal{A}$ in the Frobenius norm. 
Here the expectation is taken over the Gaussian sketches $S$, and the integers $\ell_i$ can each take any value such that 
$1<\ell_i\leq \hat r_i-r_i$. 
\end{theorem}

{\sc Proof.}
Without loss of generality assume $p_i=i$ for all $i$. Consider the first step $i=1$. 
Recall
the LS problem 
$\| \mathcal{B}^{\rm{new}} \times_i \hat F_i - \mathcal{B}^{\rm{old}} \|_F$ is equivalent to low-rank approximation of the unfolding $A_{(i)}$, and that with the sketched solution $\hat F_i$ for \eqref{eq:LSsketch}, the approximation 
$A_{(i)} \approx \hat F_i \Omega_i A_{(i)} =A_{(i)}S^T(\Omega_i A_{(i)}S^T)^\dagger \Omega_i A_{(i)}$
is the GN approximation as discussed above. Moreover, by assumption there exists a rank-$r_1$ matrix $B$ such that $\|A_{(i)}-B\|_F\leq \|\mathcal{A}-\mathcal{A}_{{\rm opt}}\|_F$. Hence by~\cite{tropp2017practical,Nakatsukasa20fast}, 
\eqref{eq:LSsketch}
is solved such that the computed $\hat F_i$ satisfies
\begin{equation}
  \label{eq:thmproof}
\mathbb{E}[\| \mathcal{B}^{\rm{new}} \times_i \hat F_i - \mathcal{B}^{\rm{old}} \|_F]\leq 
\sqrt{1+\frac{s_i}{s_i-\hat r_i-1}}\sqrt{1+\frac{\hat r_i}{\hat r_i-\ell_i-1}}
\|\mathcal{A}-\mathcal{A}_{{\rm opt}}\|_F.  
\end{equation}

Note that~\eqref{eq:thmproof} implies that in the second step $i=2$, 
the best rank-$r_2$ approximation of the unfolding $A_{(2)}$ of $\mathcal{B}^{\rm{new}}$ is now bounded by 
$\left(
\sqrt{1+\frac{s_i}{s_i-\hat r_i-1}}\sqrt{1+\frac{\hat r_i}{\hat r_i-\ell_i-1}}
\right)
\|\mathcal{A}-\mathcal{A}_{{\rm opt}}\|_F$. 

The result follows by repeatedly applying the above arguments for $i=1,\ldots, d$, noting that after $k$ steps, the best rank-$r_i$ approximation of the unfolding $A_{(i)}$ of $\mathcal{B}^{\rm{new}}$ is bounded in error by 
$\prod_{i=1}^k \left(
\sqrt{1+\frac{s_i}{s_i-\hat r_i-1}}\sqrt{1+\frac{\hat r_i}{\hat r_i-\ell_i-1}}\right)
\|\mathcal{A}-\mathcal{A}_{{\rm opt}}\|_F$. 
\hfill$\square$

It is worth highlighting the product in parenthesis in~\eqref{eq:theorembound}: this is in contrast to the corresponding bound for Randomized STHOSVD \cite{Minster20} only involves the sum, not the product, with respect to $j$. 
The situation is similar to~\cite[\S~5.2]{Minster20}, where the use of nonorthogonal projection prevented the analysis from using orthogonality of the residuals established in~\cite[Thm.~5.1]{Vannieuwenhoven12}. 
We expect the product in~\eqref{eq:theorembound} to tend to yield an overestimate (as the said orthogonality should hold approximately, albeit not exactly), and suspect that a more involved analysis would give a tighter bound.

\subsubsection{Solving LS via row subset selection}\label{sec:LSsolvesec}
An important aspect of the problem~\eqref{eq:LSsketch} is the large number $n_{(-1)}$ of right-hand sides; which makes it crucial that the sketching cost for these is kept low. In fact, experiments suggest that the use of an SRFT sketch often results in the computation being dominated by sketching the right-hand sides. 

To circumvent this and 
to enhance efficiency in \RTSMS\ we choose $S$ to be (instead of Gaussian) a column subselection matrix, i.e., $S$ is a column submatrix of identity. We suggest the use of \emph{subsampling}, i.e., $S$ is a row subset of identity. This way the cost of sketching 
(computing $S(A_{(1)}^T\Omega_1^T)$ and $SA_{(1)}^T$)
becomes minimal.

To choose the subsampled indices we use the leverage scores~\cite{drineas2012fast}, which is a common technique in randomized LS problems and beyond~\cite{murray2023randomized}. These are the squared row-norms of the orthogonal factor of the orthonormal column space of the coefficient matrix $(\Omega_1 A_{(1)})^T$, and can be approximated using sketching~\cite{drineas2012fast,mahoney2011randomized} with $O(N\hat r_1\log N)$ operations (with an SRFT sketch) for an $\hat r_1\times N$ coefficient matrix; here $N=n_{(-1)}$. 

In brief, approximate leverage scores are computed as follows~\cite{drineas2012fast}: 
First sketch the matrix to compute $Y^T(\Omega_1 A_{(1)})^T$ and its thin QR $Y^T(\Omega_1 A_{(1)})^T=QR$, where $Y$ is an $N\times O(\hat r_1)$ SRFT sketch. Then $(\Omega_1 A_{(1)})^TR^{-1}$ is well-conditioned, so we estimate its row norms via sampling, i.e., the row norms of $(\Omega_1 A_{(1)})^TR^{-1}G$ where $G$ is a Gaussian matrix with $O(1)$ columns. Importantly, the SRFT sketch $Y^T(\Omega_1 A_{(1)})^T$ (which requires $O(Nr_1\log N)$ operations) is needed also in the rank estimation process, so this computation incurs no additional cost; and note that it is cheaper than computing $A_{(1)} Y$ (i.e., sketching the right-hand side of~\eqref{eq:LSsketch} with $Y$), because $\hat r_1<n_1$. We then choose $s = O(\hat r_1\log \hat r_1)$ indices from $\{1,\ldots, z_1\}$, by randomly sampling \emph{without} replacements, with the $i$th row chosen with probability $\ell_i/(\sum_{j=1}^{z_1}\ell_j)$. We then form the resulting subsample matrix $S \in\mathbb{R}^{s\times z_1}$, which is a column-submatrix of identity $I_{n_{(-1)}}$. We then solve the subsampled problem \eqref{eq:LSsketch}. 

To understand the role of leverage scores, let us state a result on the residual for a sketched least-squares problem
~\eqref{eq:LSstate} for a general $S$ (not necessarily a subsampling matrix), which we state in terms of a standard least-squares problem $\min_X\|AX-B\|_F$.

\begin{proposition}\label{prop}
Consider the $m\times n$ LS problem $\min_X\|AX-B\|_F$ with $A\in\mathbb{R}^{m \times n} (m\geq n), B\in\mathbb{R}^{m \times n_1}$. Let $A=QR$ be the thin QR factorization with $Q\in\mathbb{R}^{m\times n}$, and consider $SQ\in\mathbb{R}^{s\times n}$ with $s\geq n$. Let $\hat X_{*}$ denote the solution for $\min_X\|S(AX-B)\|_F$. 
Then we have 
\begin{equation}  \label{eq:LSthm}
\|A\hat X_*-B\|_F\leq \frac{\|S\|_2}{\sigma_{\min}(SQ)}\min_{X}\|AX-B\|_F.
\end{equation}
\end{proposition}
{\sc Proof.}
This can be seen as a repeated application of a bound for a subsampled least-squares problem, e.g. in~\cite{chaturantabut2010nonlinear}, slightly generalized to multiple right-hand sides: 
\ignore{consider the $i$th row of~\eqref{eq:LSsubsample}, which for simplicity we write 
(in transposed form) $\min_{u\in\mathbb{R}^{\hat r_i}}\|S(Bu-b)\|_2$; here $B=(\Omega_1A_{(1)})^T$ and $b$ is the $i$th column of $A_{(1)}$. Then $\|S(Bu-b)\|_2=\|(SQ)Ru-Sb\|_2$, for which the solution is $u_*=R^{-1}(SQ)^\dagger Sb$, and hence 
\begin{align*}
\|Bu_*-b\|_2&=\|QRR^{-1}(SQ)^\dagger Sb-b\|_2=\|(I-Q(SQ)^\dagger S)b\|_2. 
\end{align*}}

Consider the $i$th column, which for simplicity we write 
$\min_{x}\|Ax-b\|_2$. Then $\|S(Ax-b)\|_2=\|(SQ)Rx-Sb\|_2$, for which the solution is $x_*=R^{-1}(SQ)^\dagger Sb$, and hence 
\begin{align*}
\|Ax_*-b\|_2&=\|QRR^{-1}(SQ)^\dagger Sb-b\|_2=\|(I-Q(SQ)^\dagger S)b\|_2. 
\end{align*}
Now note that $Q(SQ)^\dagger S$ is an (oblique) projection matrix 
$(Q(SQ)^\dagger S)^2=Q(SQ)^\dagger S$ 
onto the span of $Q$, and so $I-Q(SQ)^\dagger S$ is also a projection. It hence follows that
\begin{align*}
\|(I-Q(SQ)^\dagger S)b\|_2&=
\|(I-Q(SQ)^\dagger S)Q_\perp Q_\perp^Tb\|_2\\
&\leq \|(I-Q(SQ)^\dagger S)\|_2 \|Q_\perp^Tb\|_2
=\|Q(SQ)^\dagger S\|_2 \|Q_\perp^Tb\|_2  , 
\end{align*}
 where the last equality holds because $Q(SQ)^\dagger S$ is a projection~\cite{szyld2006many}. 

Finally, noting that 
$\|Q_\perp b\|_2=\min_{x}\|Ax-b\|_2$, 
we conclude that 
\[
\|Ax_*-b\|_2
\leq \|Q(SQ)^\dagger S\|_2 \min_{x}\|Ax-b\|_2
\leq \frac{ \|S\|_2}{\sigma_{\min}(SQ)} \min_{x}\|Ax-b\|_2. 
 \]
The claim follows by repeating the argument for every column $i=1,\ldots, n$. 
\hfill$\square$

Let us discuss Proposition~\ref{prop} when $S$ is a subsampling matrix generated via approximate leverage scores. 
First, we note that the use of leverage scores for the LS problem here is different from classical ones~\cite{woodruff2014sketching} in two ways: First, usually, leverage scores are computed for the subspace of the augmented matrix $[A,B]$, including the right-hand side (and usually there is a single right-hand side). We avoid this because this necessitates sketching the right-hand sides, of which there are many ($n_1$)
of them; 
this becomes the computational bottleneck, which is precisely why we opted to finding an appropriate row subsample to reduce the cost. 
The by-product is that the suboptimality of the computed solution is governed by $\frac{\|S\|_2}{\sigma_{\min}(SQ)}$, instead of the subspace embedding constant (which can be $<1$) as in standard methods. 
This leads to the second difference: 
 we do not scale the entries of $S$ inverse-proportionally with the leverage scores $\ell_i$, so as to avoid $\|S\|_2\gg 1$, which can happen when a row with low leverage score happens to be chosen. 
For the same reason we sample rows without replacements. 
We thus ensure $\|S\|_2=1$, which 
guarantees a good solution as long as $1/\sigma_{\min}(SQ)$ is not large\footnote{
We should be content with $\sigma_{\min}(SQ)= O(1/\sqrt{z_1})$, which is what we expect if $Q$ was Haar distributed; 
It is important to note that Proposition~\ref{prop} is an upper bound, and typically an overestimate by a factor $\approx \sqrt{z_i}$. }. 
Note that this means the standard theory for leverage scores do not hold directly; however, by choosing large rows of $Q$ with high probability we tend to increase the singular values of $SQ$, and with a modest number of oversampling we typically have a modest $(\sigma_{\min}(SQ))^{-1}$. 
If it is desirable to ensure this condition, we can use the estimate $\sigma_{i}(SQ)\approx \sigma_i(S(\Omega_1 A_{(1)})^TR^{-1})$, which follows from the fact that $(\Omega_1 A_{(1)})^TR^{-1}$ is well-conditioned with high probability, by the construction of $R$ (this is also known as whitening~\cite{nakatsukasa_tropp}). 

We note that LS problems with many right-hand sides were studied by Clarkson and Woodruff~\cite{clarkson2017low}, who show that leverage-score sampling based on $A$ (and not $b$) gives a solution with good residual. However, the failure probability decays only algebraically in $s$, not exponentially. By contrast, in standard leverage score theory, $SQ$ becomes well-conditioned with failure probability decaying exponentially in $s$~\cite[Ch.~6]{murray2023randomized}. We prefer to keep the failure probability exponentially low, given that LS has a large number of right-hand sides.

It is worth noting that leverage score sampling is just one of many methods available for column/row subset selection. Other methods include pivoted LU and QR~\cite{gu1996efficient,dong2021simpler}, and the BSS method originally developed for graph sparsification~\cite{batson2009twice}. We chose leverage score sampling to avoid the $z_in_i^2$ cost required by deterministic methods, and because Proposition~\ref{prop} and experiments suggest that oversampling (selecting more than $\hat r_i$ rows) can help improve the accuracy. 

The complexity is $O(sr_1n_{(-1)}+sr_1^2)$; $sr_1^2$ for the QR factorization $(\Omega_1 A_{(1)}S)^T=QR \in\mathbb{R}^{s\times r_1}$, and\footnote{While 
$\Omega_1 A_{(1)}S$ is simply an extraction of the columns of $\Omega_1 A_{(1)}$ specified by $S$, this step is not always negligible in an actual execution. It is nonetheless significantly faster than sketching $\Omega_1 A_{(1)}$ from the right.}
$sr_1n_1$ for computing $R^{-1}Q^T A_{(1)}$.
Together with the cost for leverage score estimation, the overall cost for solving~\eqref{eq:LSstate} is $O(n_{(-1)}\hat r_1\log n_{(-1)}+sr_1n_{1}+sr_1^2)$.

Our experiments suggest that the above process of solving~\eqref{eq:LSstate} via~\eqref{eq:LSsketch}, which can be seen as an instance of a sketch-and-solve approach for LS problems, does not always yield satisfactory results: the solution accuracy was worse by a few digits than existing algorithms such as STHOSVD. 
The likely reason is numerical instability; qualitatively, the matrix $\Omega_1 A_{(1)}$ is highly ill-conditioned, and hence the computation of its sketch also comes with potentially large relative error (that is, significantly larger than the error with $A_{(1)}QQ^T$ in~\eqref{eq:HMTapprox}; the cause is likely numerical errors, as the theory shows the residual of sketch-and-solve methods is within a modest constant of optimal~\cite{woodruff2014sketching}.). The situation does not improve with a sketch-to-precondition method~\cite{meier2023sketchandprecondition}.

In \RTSMS\ we employ  two techniques to remedy the instability: \emph{regularization} and \emph{iterative refinement}. 

\subsubsection{Regularization and refinement}
To improve the solution quality of the subsampled LS problem  
$\min_{\hat F_1}\|
S(A_{(1)}^T\Omega_1^T \hat F_1^T - A_{(1)}^T)
\|_F$, we introduce a common technique of Tikhonov regularization~\cite{hansen1998rank}, also known as ridge regression. That is, instead of~\eqref{eq:LSstate} we solve for a fixed $\lambda>0$
\begin{equation}  \label{eq:regLS}
\min_{\hat F_1}\|
S_1(A_{(1)}^T\Omega_1^T \hat F_1^T - A_{(1)}^T)\|_F^2 + \lambda \|\hat F_1\|_F^2. 
\end{equation}
This is still equivalent to an LS problem with multiple independent right-hand sides $\min_{\hat F_1}\left\|  \begin{bmatrix}S_1A_{(1)}^T\Omega_1^T  \\    \sqrt{\lambda} I  \end{bmatrix}\hat F_1^T-
  \begin{bmatrix}
    S_1 A_{(1)}^T\\0
  \end{bmatrix}
\right\|_F^2$, 
and can be solved in $O(sr_1n_{(-1)}+sr_1^2)$ operations. 
We take $\lambda=O(u\|\Omega_1 A_{(1)} S_1\|)$. Regularization is known to attenuate the effect of solution blowing up due to the presence of excessively small singular values in the coefficient matrix (some of which can be numerical artifacts). 

The second technique we employ is iterative refinement~\cite[Ch.~12]{Higham:2002:ASNA}. The idea is to simply solve the problem twice, but we found that resampling can be helpful: denoting by $\hat F_1^{(1)}$ the computed solution of~\eqref{eq:regLS}, 
we compute the residual in a \emph{different} set of subsampled columns, also obtained by the leverage scores (the two sets $S_1,S_2$ are allowed to overlap, but they differ substantially), update the (unsketched) right-hand side matrix 
$B:=A_{(1)}^T - A_{(1)}^T \Omega_1^T (\hat F_1^{(1)})^T $
and solve 
\begin{equation}  \label{eq:regLS2}
\min_{\hat F_1^{(2)}} \|S_2(A_{(1)}^T\Omega_1^T (\hat F_1^{(2)})^T - B)\|_F^2 + \lambda \|\hat F_1\|_F^2. 
\end{equation}
We then take the overall solution to be $\hat F_1^{(1)}+\hat F_1^{(2)}$. It is important in practice that the same $\lambda$ is used in~\eqref{eq:regLS} and \eqref{eq:regLS2}, even though the right-hand sides $A_{(1)}$ and $B$ are typically vastly different in norm $\|A_{(1)}\|_F\gg \|B\|_F$. This is because the goal of the second problem~\eqref{eq:regLS2} is to add a correction term, and not to solve \eqref{eq:regLS2} itself accurately. 

\begin{algorithm}[!h] 
\caption{Algorithm for the least-squares problem
$\min_{\hat F_i\in\mathbb{R}^{n_i \times \hat r_{i}}}\|(A_{(i)}^T\Omega_i^T) \hat F_i^T - A_{(i)}^T\|_F
$ as in~\eqref{eq:LSstate} arising in RTSMS. 
\\
 Inputs are $A_{(i)} \in \mathbb{R}^{z_i\times n_i}$ and  $\Omega_i\in\mathbb{R}^{\hat r_i\times z_i}$, with 
$z_i := (\Pi_{j=1}^{i-1} \hat r_j) (\Pi_{j=i+1}^{d} n_j).$
}
\label{RTSMS:LSalg}
\begin{algorithmic}[1]
     \STATE  Compute approximate leverage scores $\ell_i$: \\ \hskip2em
Using QR factorization $(\Omega_i M Y_i)^T=QR$ from Algorithm~\ref{RTSMS:alg}, 
$\ell_i$ is the $i$th row-norm of $(\Omega_i A_{(i)} Y_i)^T(R_i^{-1}G)$, for a standard Gaussian $G\in\mathbb{R}^{\hat r_i \times 5}$.
     \STATE 
If $i=p_1$, choose $s_i = 4 k \hat r_i$ indices, otherwise choose $s_i = 3k \hat r_i$ indices from $\{1,\ldots, z_i\}$, the $i$th row chosen with probaility $\ell_i^2/(\sum_{j=1}^{z_i}\ell_j^2)$ without repetition. Form the resulting subsample matrices $S_1, S_2\in\mathbb{R}^{s_i \times z_i}$.
\STATE  Use Tikhonov regularization with $S_1$ to compute $\hat F_i^{(1)}$ solving \eqnref{eq:regLS}.
     \STATE  Use iterative refinement with $S_2$ to compute $\hat F_i^{(2)}$ solving \eqnref{eq:regLS2}. 
     \STATE Compute $\hat F_i := \hat F_i^{(1)} + \hat F_i^{(2)}$.
\end{algorithmic}
\end{algorithm}

It is perhaps surprising that a standard sketch-and-solve (even with regularization) does not always yield satisfactory solutions. While the algorithm presented here always gave good computed outputs in our experiments, it is an open problem to prove its stability, or the lack of it (in which case alternative methods that would guarantee stability). The stability of randomized least-squares solvers is generally a delicate topic~\cite{meier2023sketchandprecondition}, particularly when the coefficient matrix is numerically rank deficient.

\subsection{HOSVD variant of RTSMS}\label{RHOSVDSMS:subsec}
Our algorithm can be complemented with a further conversion to HOSVD such that factor matrices have orthonormal columns and the core tensor has the so-called all-orthogonality property~\cite{de2000multilinear}. One way to do this standard deterministic step is to apply Alg.~\ref{Tucker2HOSVD:alg} in supplementary materials which has the flexibility of choosing whether a multilinear singular value thresholding should also be carried out. In our numerical experiments we denote this variant with \texttt{RHOSVDSMS}.

\subsection{Fixed-rank variant of RTSMS}
A significant aspect of RTSMS is its rank adaptivity; it can automatically adjust the numerical multilinear rank of the tensor according to a given tolerance for the relative residual. 
Here we discuss a variant of RTSMS in which the multilinear rank is assumed to be known a priori. Most of the algorithms for the Tucker decomposition in the literature are of this type. In our numerical experiments we denote this variant with \texttt{\RTSMS-FixedRank} which is essentially Steps 20--23 of RTSMS. If desired, this variant can also be converted to the HOSVD form without thresholding.

\subsection{Computational complexity}
In Table~\ref{costs:tab} we summarize the number of arithmetic operations involved in algorithms for computing a Tucker decomposition. For simplicity we assume that the order-d tensor $\mathcal{A}$ is $n \times n \dots \times n$ and the target rank is $r \times r \dots \times r$ (we do not include the cost for rank estimation; it is usually comparable to the algorithm itself). 
In addition, we assume without loss of generality that the processing order is $1,2,\dots, d$. The cost of the single-pass Tucker algorithm \cite{Sun20} is obtained by taking the factor sketching parameters to be $k = 3/2r$ and the core sketching parameters to be $s = 2k = 3r$. 

\small
\begin{table}[!h]
\begin{center}
\caption{Computational complexity of fixed-rank algorithms for computing rank $(r, r, \dots, r)$ Tucker decomposition of an order-d tensor of size $n \times n \dots \times n$, and $r\ll n$. 
$\hat r = r+p$ where $p$ is the oversampling factor, e.g. $p=5$ or $p=0.5r$. \label{costs:tab}}
\begin{tabular}{l|cll}
algorithm & dominant  & sketch  & dominant operation  \\
 & cost  & size & \\
\hline HOSVD & $d n^{d+1}$ & &\small SVD of d unfoldings each of size \\
  \cite{Tucker66, de2000multilinear} & & &\small $n \times n^{d-1}$ \\
\hline STHOSVD  & $n^{d+1}$  & &SVD of $A_{(1)}$ which is $n \times n^{d-1}$. \\
\cite{Vannieuwenhoven12}  & & & (Later unfoldings are smaller\\
  & & & due to truncation)\\
\hline R-HOSVD & $d r n^d$ & 
$\hat r\times n^{d-1}$
& 
computing  $A_{(i)} \Omega_i$ where $\Omega_i$ of \\
\cite{Minster20}
 &   \cite[Tab. 1]{Minster20} & & size $n^{d-1} \times \hat r$  and then forming \\
&  & & $Q_i^T A_{(i)}$  for all $i$  \\
\hline R-STHOSVD & $r n^d$ & $\hat r\times n^{d-1}$ & forming $A_{(1)} \Omega_1$ with $\Omega_1$ of size \\
\cite{Minster20,Zhou14}  &   
\cite[Tab. 1]{Minster20}&   & $n^{d-1} \times \hat r$. Subsequent unfoldings  \\ 
 & & & and sketching matrices are smaller \\  
\hline single-pass &$r n^d$  & $\hat r\times n^{d-1}$ & sketching by structured (Khatri   \\
\cite{Sun20} (also \cite{Malik18}) & \cite[Tab. 2]{Sun20} & & -Rao product) dimension   \\
 &   &  & reduction maps\\
\hline RTSMS & $r n^d$  & $\hat r\times n$ & computing $\Omega_1 A_{(1)}$ with $\Omega_1$ of size  \\ 
& ($n^d\log n$)  & &  $\hat r  \times n^{d-1}$  \end{tabular}
\end{center}
\end{table}
\normalsize 

It is evident that randomization reduces the exponent of the highest order term by one, and sequential truncation reduces the corresponding coefficient. 

While Table~\ref{costs:tab} does not immediately reveal a cost advantage of \RTSMS\ over R-STHOSVD and single-pass Tucker, let us repeat that the single-mode nature of the sketching can be a significant strength. For example, when the sketches are taken to be Gaussian, the cost of generating the sketches is lower with \RTSMS\ by a factor $O(n^{d-2})$. 
Moreover, while we mainly treat Gaussian sketches, one can reduce the complexity by using structured sketches; e.g. with an SRFT sketch the complexity becomes $O(n^{d}\log n)$ as listed in parenthesis\footnote{In theory~\cite[\S.~9.3]{MartinssonTroppacta} this can be 
reduced to $O(n^{d}\log r)$, resulting in strictly lower complexity than $O(n^{d}r)$. However, the corresponding implementation of the fast transform is intricate and often not available. One can also use sparse sketches~\cite{woodruff2014sketching} to get $O(n^{d})$ complexity.} in the table, which can be lower than $rn^d$. Such reduction is not possible with other methods based on finding the (orthonormal) factor matrices first, 
because the computation of $\mathcal{A}\times_1 F_1$ necessarily requires $rn^d$ operations, since $F_1$ is generally dense and unstructured. Another advantage of \RTSMS\ is that computing $\Omega_1\mathcal{A}_{(1)}$ requires far less communication than $\mathcal{A}_{(1)}\Omega_2$, as with $\Omega_1\mathcal{A}_{(1)}$ the local computation (where $\mathcal{A}_{(1)}$ is split columnwise) is direcly part of the output.

The table shows a somewhat simplified complexity; for example, \RTSMS\ also requires $O(n^{d-1}rd\log n)$ for finding the leverage scores, which is usually no larger than $n^{d}r$ but can be dominant when $d$ is large. 

\subsection{Implementation details}
Let us address a few points concerning the specifics of our RTSMS implementation\footnote{Our MATLAB implementation is available at \url{https://github.com/bhashemi/rtsms}.}. 
\begin{itemize}

\item We noticed that the final residual depends more strongly on the error made in the least-squares problem in computing the first factor matrix as opposed to the later modes. An analogous observation was made by Vannieuwenhoven, Vandebril and Meerbergen \cite[sec. 6.1]{Vannieuwenhoven12}, who emphasized the influence of the error caused by the first projection on the quality of the STHOSVD approximation in comparison with the remaining projections. For this reason we take a slightly larger number of oversamples when processing the first mode.

\item Tensor-matrix contraction is a fundamental operation in tensor computation. Such mathematical operations involve tensor unfoldings and permutations which are memory-intensive. It is therefore advantageous to decrease the need for unfolding and permutation of large tensors so as to reduce data communication through the memory hierarchy and among processors. Conventional implementations of tensor-matrix contractions involve permutation of the tensor (except for mode-1 contraction) so that the computation can be performed by calling BLAS; see the function \texttt{tmprod} in TensorLab, for instance. While our implementations employ TensorLab, we adapt the use of \texttt{tensorprod}, introduced in MATLAB R2022a, in order to accelerate tensor-matrix contractions by avoiding explicit data permutations.

\end{itemize}

\paragraph{Further rank truncation}We explored the following two strategies to further truncate the computed Tucker decomposition. 
\begin{itemize}
\item The first strategy can be applied {\em within} Step 22
of Alg.~\ref{RTSMS:alg}. 
Namely, if it is found that 
$\| \mathcal{B}^{\rm{new}} \times_i \hat F_i - \mathcal{B}^{\rm{old}} \|_F$ can be reduced below the tolerance with a rank lower than $\hat r_i$, one can do so with negligible extra cost.
Once it is found that $\hat F_i$ can be of rank $\ell<\hat r$ so that $\hat F_i= F_iG_i$ where $F_i$ has $\ell$ columns, we accordingly reduce $(A_{(1)}^T\Omega_1^T) \hat F_i^T$ to 
$(A_{(1)}^T\Omega_1^TG_i^T) F_i^T$. This results in a truncation of the current core tensor in the $i$-th mode.\item The second strategy is based on the HOSVD and is applied once the execution of Alg.~\ref{RTSMS:alg} is completed.
 We first use QR factorizations of the factors $F_i$ and deterministic STHOSVD to convert the Tucker decomposition 
\[
\mathcal{A} \approx \llbracket \mathcal{C}; F_1 , F_2 ,\dots, F_d \rrbracket
\]
to a HOSVD of the form
\[
\mathcal{A} \approx  \llbracket \mathcal{\check C}; U_1 , U_2 ,\dots, U_d \rrbracket.
\]
We then compute the higher order modal singular values of $\mathcal{A}$ from $\mathcal{\check C}$ and then compare those modal singular values with the input tolerance to decide where to truncate the factors as well as the core tensor. See Alg.~\ref{Tucker2HOSVD:alg} in supplementary materials.
\end{itemize}

The second strategy brings the Tucker decomposition into HOSVD and applies a multilinear singular value thresholding operator, commonly employed in the context of low-rank matrix recovery; see~\cite{singValThresh} for instance. 
In our numerical experiments we observed that this second strategy gives results whose multilinear rank and accuracy have better consistency with changes in the given tolerance compared with the first strategy. Details are available in a pseudocode in supplementary materials (Alg.~\ref{Tucker2HOSVD:alg}).

\section{Experiments}\label{exper:sec}

In all the following experiments we use parameters detailed here. The oversampling parameter is set to be $\tilde p = 5$ in all randomized techniques of \cite{Minster20} and \cite{Malik18}. In the randomized techniques of \cite{Minster20}, we use $[1,2, \dots, d]$ as the processing order vector of modes In the Tucker-TensorSketch \cite{Malik18}, we set the sketch dimension parameter $K$ to be equal to $\tilde p$ above, so $K = 5$, a tolerance of $1\times 10^{-15}$, and the maximum number of iterations to be $50$. In our experiments we always take the sketches to be Gaussian, as the dimensions $n_i$ and especially $r_i$ are not large enough for other sketches (e.g. SRFT) to outperform it. We set the oversampling parameter $k=4$. All experiments were carried out in MATLAB 2023a on a machine with 512GB memory.

We repeat each experiment for computing decomposition of the form $\mathcal{\tilde A} =  \llbracket \mathcal{C}; U_1 , U_2 ,\dots, U_d \rrbracket$ five times and report the average CPU time (the variance is not large) as well as the geometric mean of the relative residuals defined by 
\[
 {\rm relative\ residual} = \frac{\| \mathcal{A} - \mathcal{\tilde A} \|_F}{\|\mathcal{A}\|_F}\cdot
\]
As our main algorithm is adaptive in rank, in our first four examples we focus on experiments in which an input tolerance is specified and, in addition to the average CPU time and residual, we also report the average numerical multilinear rank (rounded to the nearest integer) as computed by each algorithm. In those experiments we compare RTSMS and RHOSVDSMS (see section~\ref{RHOSVDSMS:subsec}) with adaptive-rank R-STHOSVD (Alg.~\ref{adapt_RST:alg}), whose Step 3 incorporates \texttt{svdsketch}. Note that while the only mandatory input to \texttt{svdsketch} is an $m \times n$ input matrix $A$, it allows specifying the following input parameters:
The input tolerance tol,  maximum subspace dimension, blocksize, maximum number of iterations, and number of power iterations performed (default value: 1). Among these, it is worth noting that the input tolerance tol is required to satisfy
  \begin{equation}    \label{eq:tolrequire}
  \sqrt{{\rm machine\ epsilon}} \approx 1.5 \times 10^{-8} \leq {\rm tol} < 1,  
  \end{equation}
since \texttt{svdsketch} does not detect errors smaller than the square root of machine epsilon; see Theorem 3 and Remark 3.3 in \cite{Yu18}.  The default value of tol is ${\rm machine\ epsilon}^{1/4} \approx 1.2 \times 10^{-4}$ but note that in this paper we will specify different values of tol satisfying \eqref{eq:tolrequire}.
Also, in our tensor context the maximum subspace dimension, blocksize and the maximum number of iterations are typically $m$, $\lfloor 0.1 m\rfloor$ and 10, respectively.

\subsection{Rank-adaptive experiments}
 
\begin{example}\label{runge:ex} \normalfont
We take $\mathcal{A}$ to be samples of the Runge function $f(x,y,z) = 1/(5+x^2 + y^2 + z^2)$ on a grid of size $600 \times 600 \times 600$ consisting of Chebyshev points on $[-1, 1]^3$. Figure~\ref{runge:fig} reports the results, in which the left panel shows the speed of our algorithm over R-STHOSVD. 

\begin{figure}[!h]
\begin{center}
\includegraphics[width=.8\textwidth]{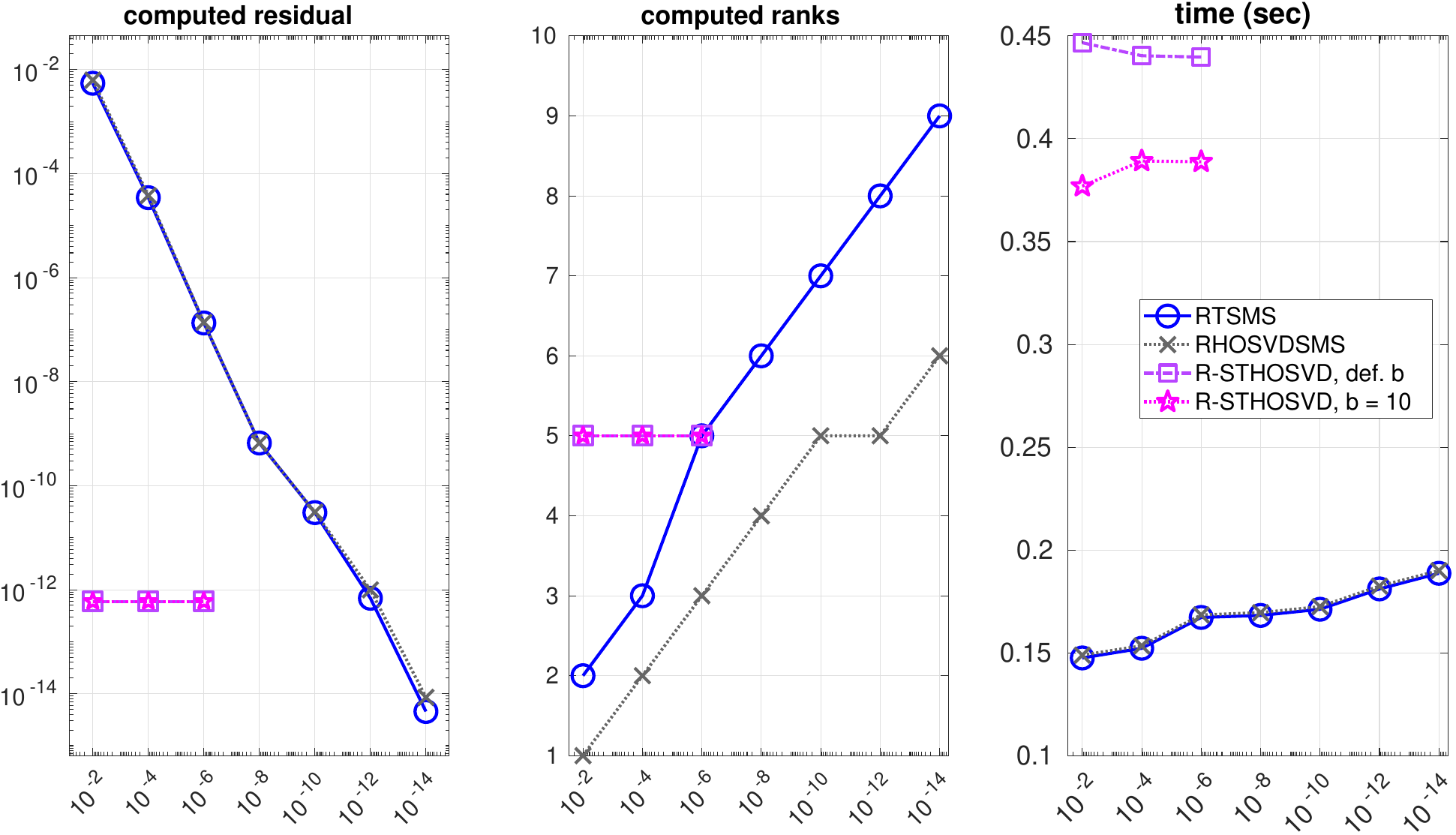}
\caption{Comparison of rank-adaptive methods in terms of the the residual (left), average multilinear ranks (middle) and time (right) for a tensor $\mathcal{A}$ of size $600 \times 600 \times 600$ in Example~\ref{runge:ex}. The horizontal axes are all the input tolerance. \label{runge:fig}}
\end{center}
\end{figure}

The default values of the block size $b$ used by \texttt{svdsketch} in R-STHOSVD are $\lfloor 0.1n_i \rfloor = 60$ in modes $i = 1$ and $2$. The corresponding average of the computed multilinear ranks is 5 (not 60 or 10 in the two executions of R-STHOSVD) reflecting the fact that orthogonalizations within randomized matrix range finders performed by MATLAB \texttt{orth} only keep $r$ columns of the unfoldings where $r$ is the computed rank. When it comes to the third mode, the default value of the block size $b$ used in R-STHOSVD is only 25, as due to its sequential truncation aspect, the third unfolding matrix turns out to be of size $n_3 \times (r_1 r_2) = 600 \times 25$. In contrast to the multilinear ranks being always equal to 5 in both executions of R-STHOSVD in this example, we observe a smooth increase in the computed ranks as the input tolerance is decreased in both RTSMS and RHOSVDSMS.

\end{example}

\begin{example}\label{wagon:ex}\normalfont
We take $\mathcal{A}$ to be samples of the Wagon function on a grid of size $800 \times 1200 \times 300$ consisting of Chebyshev points on $[-1, 1]^3$. The function appears in a challenging global minimization problem and is most complicated in the second variable which is why we took the second dimension of $\mathcal{A}$ to be the largest. See \cite{Bornemann04} for details. 

Figure~\ref{wagon:fig} illustrates our comparisons. Relative residuals in both RTSMS and RHOSVDSMS drops to about machine epsilon already for input tolerance as big as $10^{-6}$ because Wagon's function, while challenging in its three variables, is intrinsically  of low multilinear rank.

\begin{figure}[!h]
\begin{center}
\includegraphics[width=.8\textwidth]{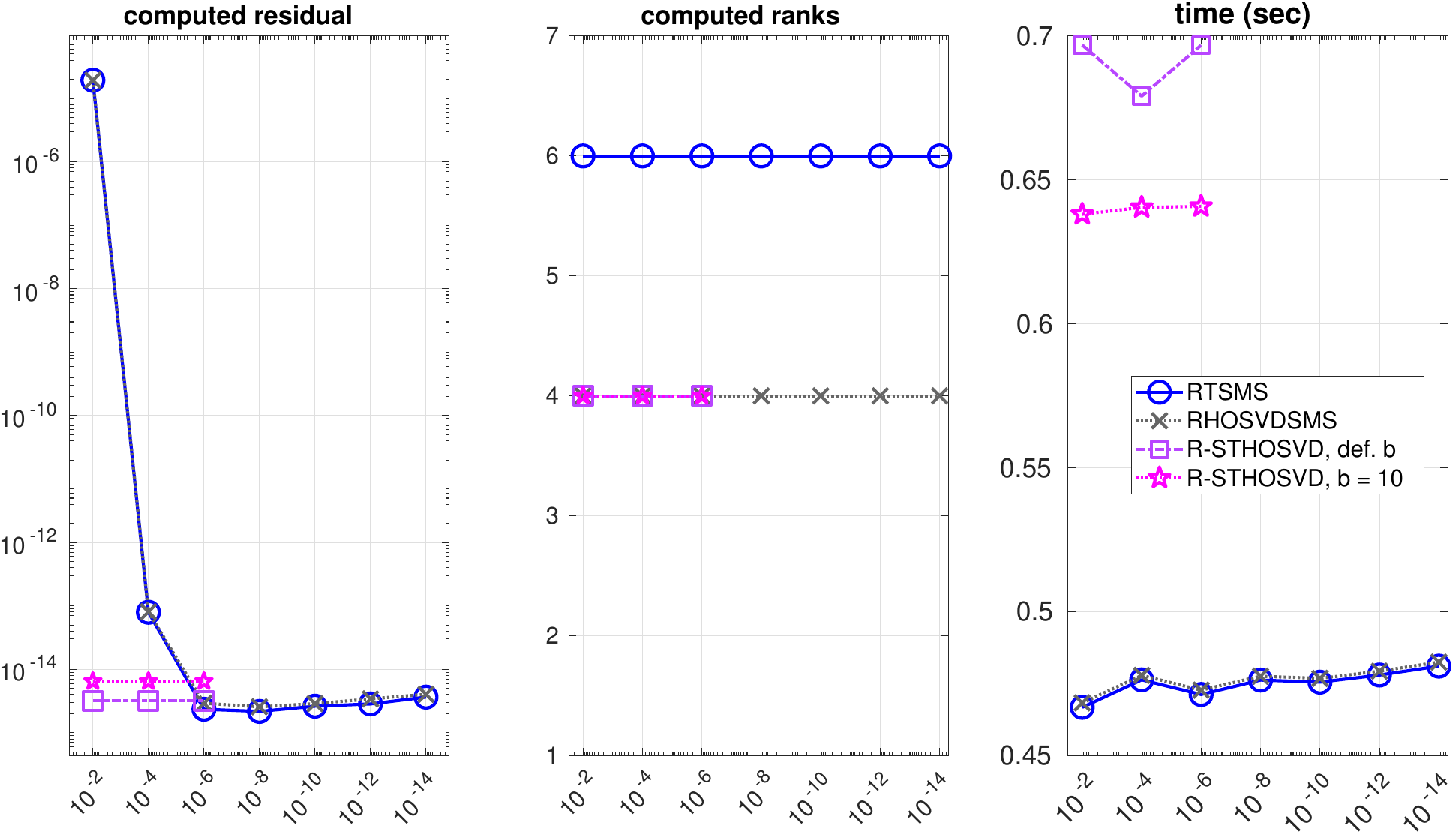}
\caption{Comparison of rank-adaptive methods in terms of the residual (left), average multilinear ranks (middle) and computing time (right) for a tensor $\mathcal{A}$ of size $800 \times 1200 \times 300$ in Example~\ref{wagon:ex}. The horizontal axes are all the input tolerance.
\label{wagon:fig}}
\end{center}
\end{figure}

\end{example}

\begin{example}\label{octant:ex}\normalfont
We take $\mathcal{A}$ to be samples of the function $f(x,y,z) = \sqrt{x^2 + y^2 + z^2}$ on a Chebyshev grid of size $1000 \times 1000 \times 1000$ in $[-1, 1]^3$. $f$ is called the Octant function in \texttt{cheb.gallery3} in Chebfun3 \cite{hashemichebfun3} and has nontrivial ranks.

See Fig.~\ref{octant:fig}. The size $b$ of blocks used by \texttt{svdsketch} in R-STHOSVD is $\lfloor 0.1n_i \rfloor = 100$ for $i=1,2,3$, and the corresponding average of the computed multilinear ranks is 22.

\begin{figure}[!h]
\begin{center}
\includegraphics[width=.8\textwidth]{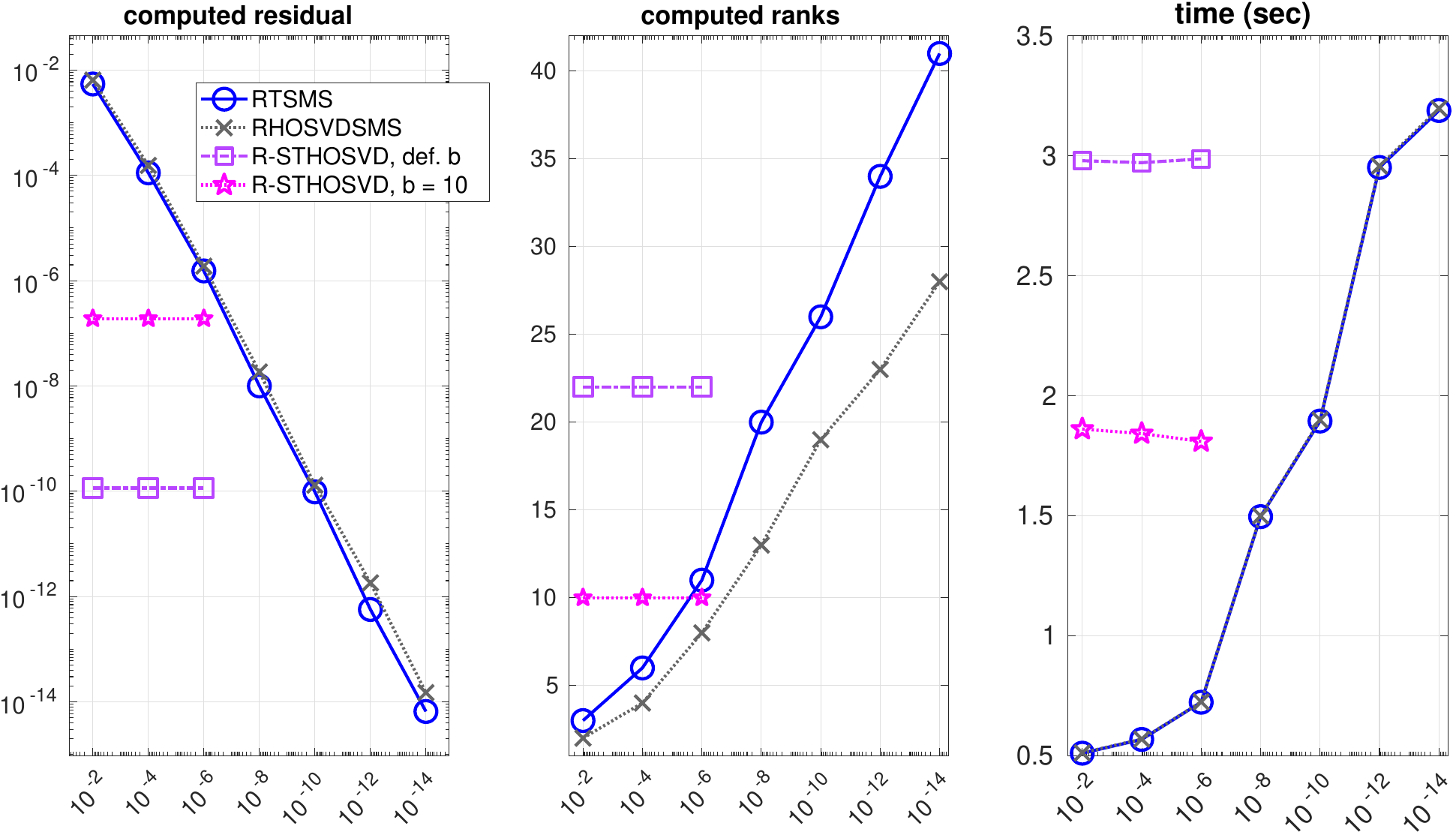}
\caption{Comparison of methods in terms of residual (left), average multilinear ranks (middle) and computing time (right) for a tensor $\mathcal{A}$ of size $1000 \times 1000 \times 1000$ in Example~\ref{octant:ex}. The horizontal axes are all the input tolerance. \label{octant:fig}}
\end{center}
\end{figure}

\end{example}

\begin{example}\label{issVidGrayscale:ex} \normalfont
We take $\mathcal{A}$ to be a 3D tensor of size $483 \times 720 \times 1280$ containing 483 frames of a video from the international space station. It corresponds to the first 16 seconds of a color video\footnote{\url{https://www.youtube.com/watch?v=aIkWx6HGol0} retrieved September 13, 2023. 
} 
which can be represented as a 4D tensor. However, we converted the color video to grayscale and only took the first 483 frames, creating an order-three tensor $\mathcal{A}$. 

We try RTSMS with two tolerances $10^{-2}$ and $10^{-3}$. In the first case, a Tucker decomposition of multilinear rank $(80, 117, 131)$ is computed in 8.7 seconds with an actual relative residual of $2.2 \times 10^{-1}$. With tol $=10^{-3}$ we get a Tucker decomposition of rank $(362, 552, 659)$ after 95.7 seconds with an actual relative residual of 
$3.9 \times 10^{-2}$. In this example we did not apply either of the two truncation strategies explained in Section~\ref{sec:RTSMS}. The variant of R-STHOSVD which, instead of the rank, takes a tolerance as input did not give an output after 5 minutes following which we stopped its execution. 

\begin{figure}[!h]
\begin{center}
\includegraphics[width=1\textwidth]{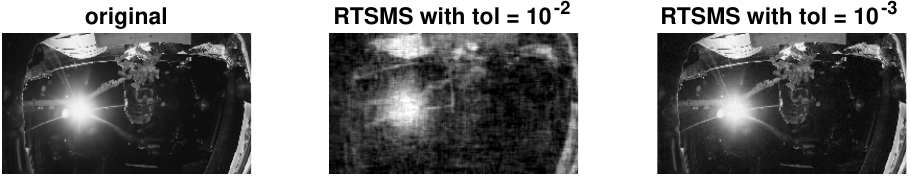}
\includegraphics[width=1\textwidth]{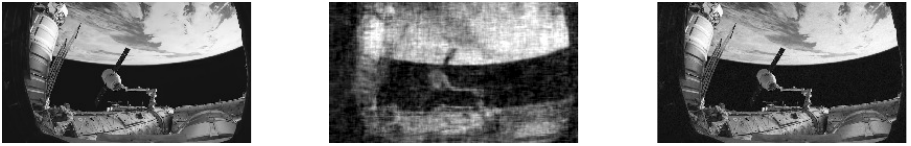}
\includegraphics[width=1\textwidth]{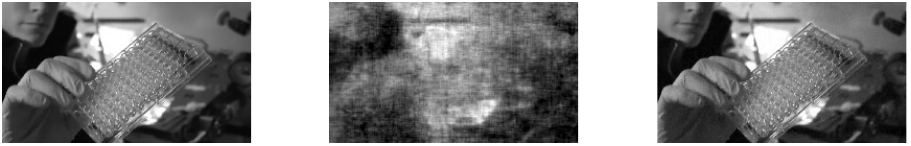}
\caption{Frames 1 (top), 200 (middle) and 400(bottom) from the International Space Station video and approximations obtained by RTSMS with two different tolerances. See Example \ref{issVidGrayscale:ex}.}
 \label{issVidGrayscale:fig}
\end{center}
\end{figure}
The images are shown in Figure~\ref{issVidGrayscale:fig}. A 16-second video is available in supplementary materials showing this comparison for all 483 frames. 

\end{example}

\begin{example}\label{miranda:ex} \normalfont
In our next example we work with a 3D tensor of size $2048 \times 256 \times 256$ from the Miranda Turbulent Flow tensor data from the Scientific Data Reduction Benchmark (SDRBench) \cite{Cabot06, Zhao20}. We acquired the dataset following the methodology outlined in \cite{Ballard22} following which we apply the HOSVD variant of RTSMS with four tolerances, specifically $10^{-2}, 10^{-3}, 10^{-4}$ and $10^{-5}$. 

We report our results in Table~\ref{miranda:tab}. See also Figures~\ref{miranda_compare:fig} and~\ref{miranda:fig} for illustrations. Here, the relative compression is the ratio of the total size of the original tensor $X$ and the total storage required for the Tucker approximation. 

\begin{figure}[!h]
\begin{center}
\includegraphics[width=0.4\textwidth]{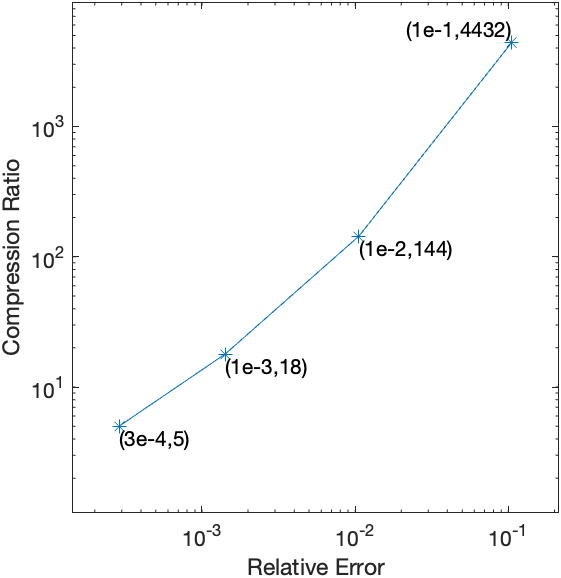}
\caption{Compression versus error obtained with HOSVDSMS with different tolerances in Example~\ref{miranda:ex}.}
\label{miranda_compare:fig}
\end{center}
\end{figure}

We obtain a compression ratio of 5X requiring $20\%$ of the size of the original tensor when the tolerance is set to $10^{-5}$. On the other hand with a tolerance of $10^{-2}$, a compression ratio of 4432 is achieved requiring only $0.02\%$ of the storage required for the original tensor. 

\begin{table}[!h]
\begin{center}
\caption{Relative error and compressions achieved for different tolerances in Example~\ref{miranda:ex}
\label{miranda:tab}}
\begin{tabular}{c|cccc}
tol & rel. error & Tucker rank & compression ratio & \% of original size \\
\hline
$10^{-2}$ & $1.05 \times 10^{-1}$ & (12,  9, 9)  & 4432 &  0.02\\
$10^{-3}$ & $1.06 \times 10^{-2}$ & (182, 54, 54)  & 144 & 0.7\\
$10^{-4}$ & $1.42 \times 10^{-3}$ & (514, 113, 110)  & 18 & 6\\ 
$10^{-5}$ & $2.88 \times 10^{-4}$ & (877, 172, 165)  & 5 & 20\\ 
\end{tabular}
\end{center}
\end{table}

\begin{figure}[!h]
\begin{center}
\includegraphics[width=0.9\textwidth]{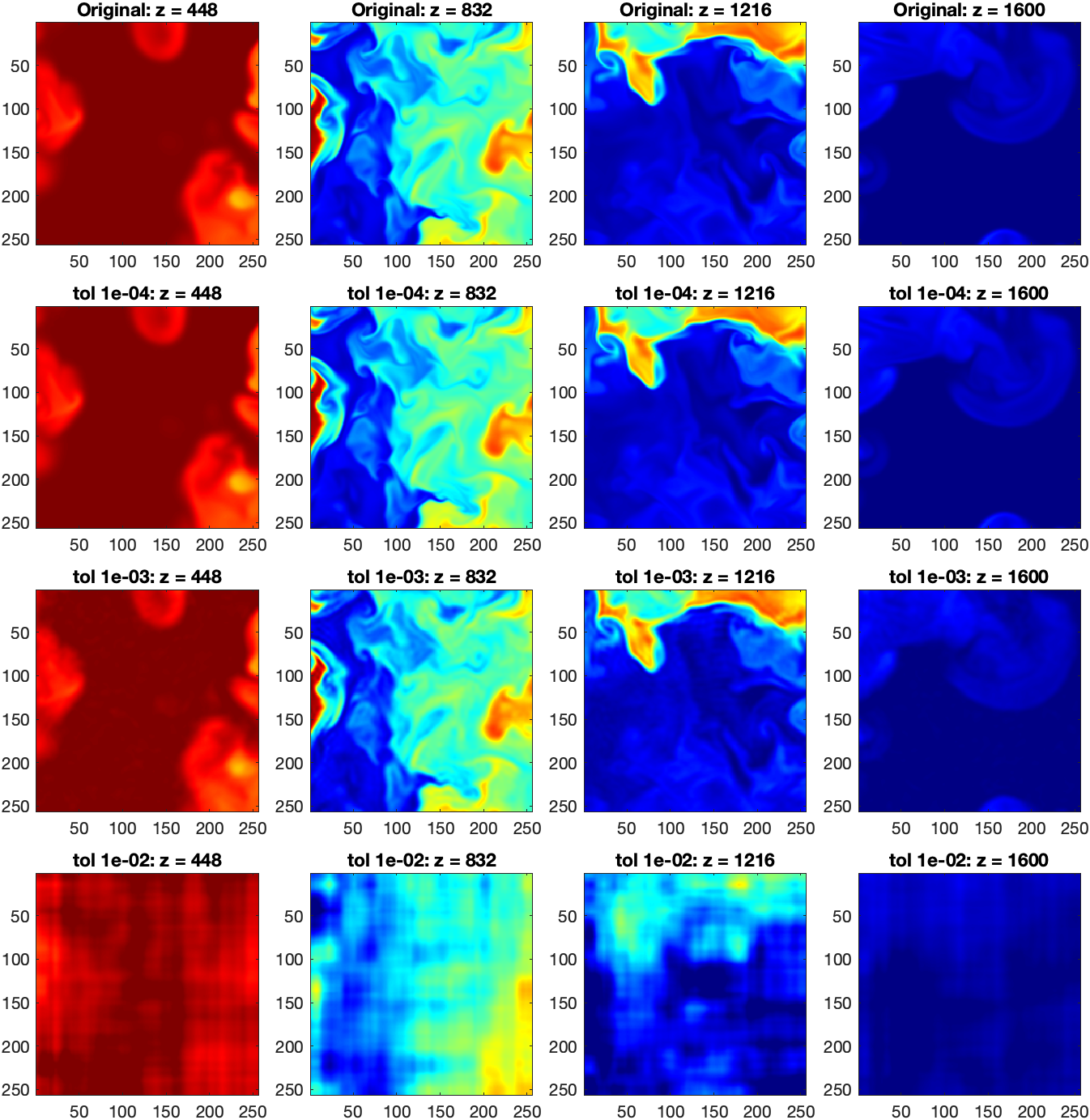}
\caption{Visualization of four slices in the xy-plane of the original density tensor and compressed representations obtained with HOSVDSMS in the numerical simulation of flows with turbulent mixing. Images produced with the tolerance of $10^{-5}$ look the same as those for $10^{-4}$ hence not displayed here. See Example~\ref{miranda:ex}.}
\label{miranda:fig}
\end{center}
\end{figure}

\end{example}

\subsection{Fixed-rank experiments}
In the following examples we examine algorithms which require the multilinear rank as input. 
While because of oversampling,  the numerical Tucker rank of the approximations computed with RTSMS is more than the input rank $\bf{r}$, RHOSVDSMS truncates those approximations to the rank $\bf{r}$. Like before we run each example five times and report the average time and residual.

\begin{example} \label{hilbert4:ex}\normalfont
We consider the four dimensional Hilbert tensor of size $150 \times 150 \times 150 \times 150$ whose entries are
\[
h_{i,j,k,l} = \frac{1}{i+j+k+l-3}\cdot
\]
We compute Tucker decompositions of multilinear rank $(r,r,r,r)$ for the following six values of $r:= 5, 10, 15, 20, 25, 30$. The results are depicted in Fig.~\ref{hilbert4:fig}. 

\begin{figure}[!h]
\begin{center}
\includegraphics[width=0.6\textwidth]{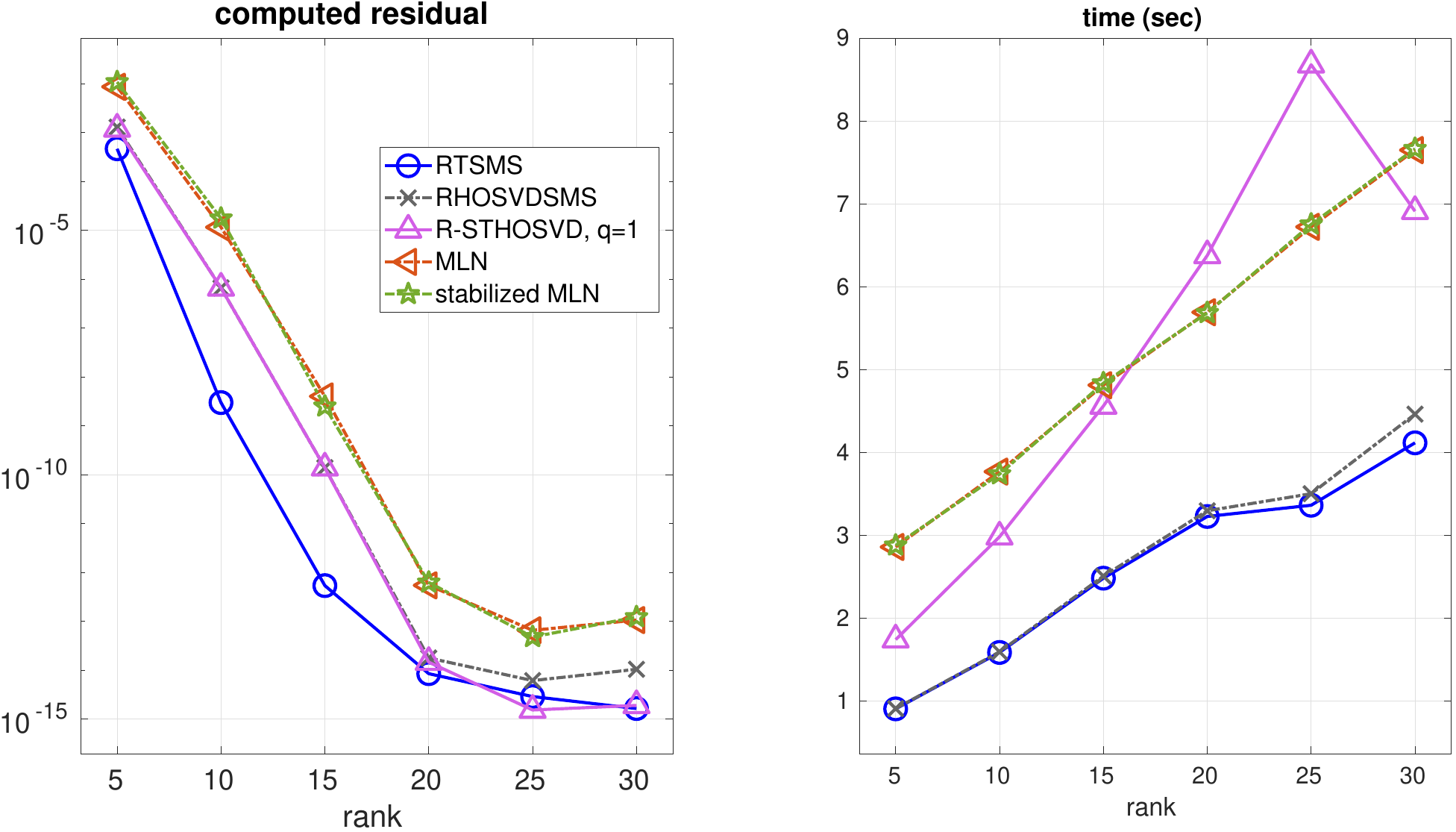}
\caption{Comparison of fixed-rank algorithms for 4D Hilbert tensor in terms of residual (left) and time (right) in Example~\ref{hilbert4:ex}.}
\label{hilbert4:fig}
\end{center}
\end{figure}

In addition to RTSMS, RHOSVDSMS, and R-STHOSVD (Alg.~\ref{Fixed_RSTHOSVD:alg}) we also plot the results obtained by the multilinear generalized Nystr\"om (MLN) and its stabilized variant \cite{Bucci23}. All the methods are comparable with the RTSMS and RHOSVDSMS slightly better both in terms of accuracy and speed. 

\end{example}

\begin{example}\label{synthetic:ex}\normalfont Motivated by demo2 in the implementations accompanying \cite{Malik18}, we construct synthetic data with noise as follows. Using Tensor Toolbox \cite{Bader06}, we generate four $n \times n \times n$ tensors of true rank $(r, r, r)$ where $n = 250, 500, 750, 1000$ and $r = 10, 12, 14, 16$, respectively. Then, a Gaussian noise at the level of $10^{-7}, 10^{-6}, 10^{-5}, 10^{-4}$ is added to the tensors, respectively. More precisely, a noise at the level of $10^{-7}$ is added to the smallest tensor ($n = 250$) and a similar noise at the level of $10^{-4}$ is applied to the largest tensor ($n = 1000$). We then call different methods to compute Tucker decomposition of noisy tensors repeating each experiment five times as before. Tucker-TS and Tucker-ALS require a few parameters which we set as follows. The tolerance and target ranks are set to the same noise level and true ranks when generating each tensor as mentioned above. In addition, as recommended in \cite{Malik18}, we set sketch dimensions to $J_1 = Kr^2$ and $J_2 = Kr^3$ with $K = 10$. 

\begin{figure}[!h]
\begin{center}
\includegraphics[width=0.6\textwidth]{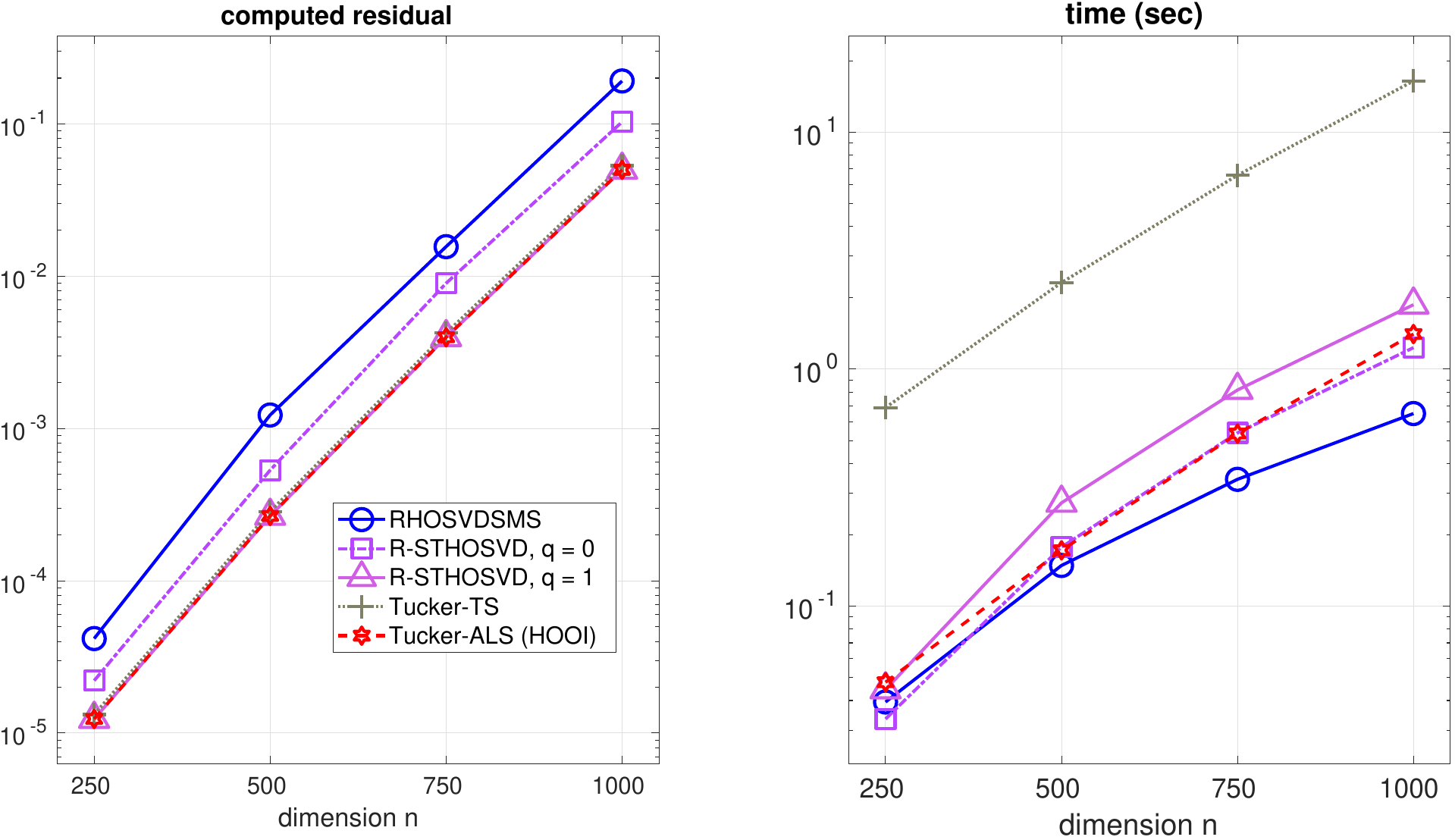}
\caption{Comparison of fixed-rank algorithms for noisy synthetic data in terms of residual (left) and time (right) in Example~\ref{synthetic:ex}. 
\label{synthetic:fig}}
\end{center}
\end{figure}

Figure~\ref{synthetic:fig} reports the outcome. While RHOSVDSMS gives a residual which is worse than R-STHOSVD by a factor of four, it is the fastest. In particular for the largest tensor, the average computing time of $110.5$ and $21.2$ seconds in Tucker-TS and R-STHOSVD, respectively 
is reduced to $4.8$ seconds. 

We also note that we encountered errors in MATLAB when computing Khatri-Rao products involved in Tucker-TS (described in section~\ref{TS:subsec}), complaining that memory required to generate the array exceeds maximum array size preference. This happens e.g., with $n = 1000$ and a rank as small as $r = 20$ in which case Tucker-TS generates an array which requires 36 GB of memory. No such errors arise with R-STHOSVD, RTSMS and RHOSVDSMS. We tried Tucker-TTMTS as well, but it gave residuals at the constant level of $10^{-1}$ for all the four tensors, and that is why Tucker-TTMTS is not shown in the plots. 

\end{example}

\begin{example}\label{medImage:ex} \normalfont
We take $\mathcal{A}$ to be the tensor of 80 snapshots from a computerized tomography (CT) kidney dataset of images\footnote{\url{https://www.kaggle.com/datasets/nazmul0087/ct-kidney-dataset-normal-cyst-tumor-and-stone}}. More specifically, $\mathcal{A}$ is a tensor of size $512 \times 512 \times 80$ corresponding to images number 328 to 407 from cyst directory in the dataset. We set the Tucker rank to be $(250, 250, 50)$. See Figure~\ref{medImage:fig} suggesting that visually R-STHOSVD with one power iteration and RTSMS give comparable approximations to the original images. In fact, RTSMS gives a mean relative residual of $1.21\times 10^{-1}$ compared with $1.07 \times 10^{-1}$ with R-STHOSVD. The average time taken by R-STHOSVD and RTSMS is 4.05 and 2.60 seconds, respectively. We see that the images are approximated with roughly the same quality by all algorithms.

\begin{figure}[!h]
\begin{center}
\includegraphics[width=.8\textwidth]{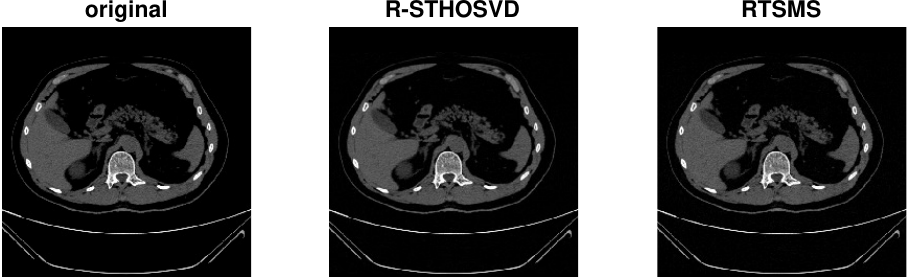}\\
\includegraphics[width=.8\textwidth]{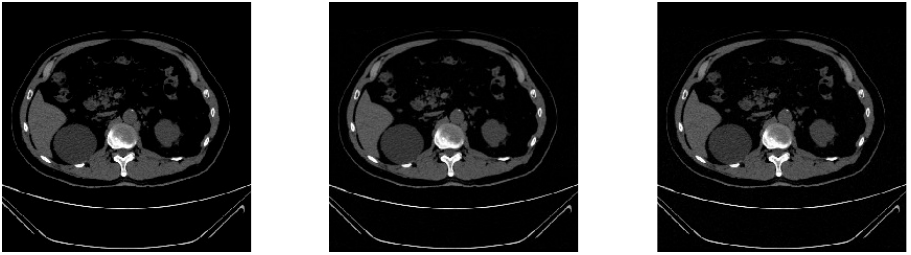}\\
\ \includegraphics[width=.8\textwidth]{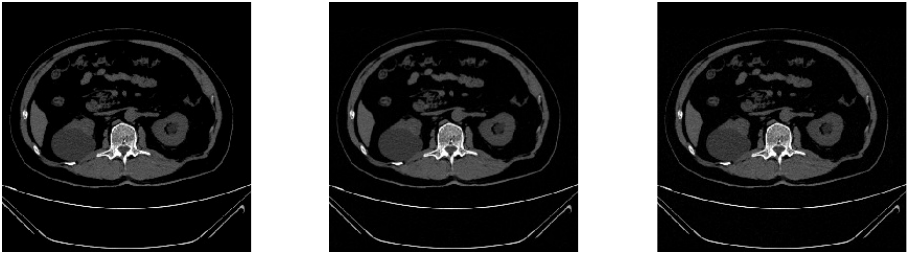}
\caption{Original images (left) and approximations obtained by R-STHOSVD (middle) and \RTSMS\ (right) corresponding to snapshots 1 (top row), 25 (middle row) and 50 (bottom row). See Example \ref{medImage:ex}. \label{medImage:fig}}
\end{center}
\end{figure}

\end{example}

\subsection*{Acknowledgments}
We would like to thank Tammy Kolda for the insightful discussions and helpful suggestions, including the experiments in Figure~\ref{miranda:fig}.

\bibliographystyle{siamplain}
\bibliography{bib2}

\newpage
\ \\
{\bf Supplementary materials.} 
Here we first recall the basic ideas of HOSVD and STHOSVD whose foundations are key to RTSMS. Also, for the sake of completeness, we present standard deterministic algorithms for HOSVD and STHOSVD as well our algorithm for converting Tucker decompositions computed by RTSMS to HOSVD. 

 Let $n_i \gg r_i$ for $i=1,2,\dots, d$ and $U_i \in \mathbb{R}^{n_i \times r_i}$ be {\em any} full-rank matrix whose columns span the column space of $A_{(i)}$ which is a subspace of $\mathbb{R}^{n_i}$. Hence, each $U_i$ has a left-inverse, i.e., $U_i^\dagger U_i = I_{r_i} \in \mathbb{R}^{r_i \times r_i}$ and  $U_i U_i^\dagger =: P_i \in \mathbb{R}^{n_i \times n_i}$ is a projection onto the column space of $U_i$ and $A_{(i)}$. Hence, $P_i U_i = U_i$ and $P_i A_{(i)} = A_{(i)}$. Using Definition 
\ref{modalProd:def}, the latter formula can be rewritten as 
\[
A_{(i)} = (\mathcal{A} \times_i P_i)_{(i)}.\]
Repeating this d times, we therefore have
\begin{align}
\mathcal{A} &= \mathcal{A} \times_1 P_1 \times_2 P_2  \dots \times_d P_d \nonumber \\
& = \mathcal{A} \times_1 U_1 U_1^\dagger \times_2 U_2 U_2^\dagger \dots \times_d U_d U_d^\dagger \nonumber\\
& = (\mathcal{A} \times_1 U_1^\dagger \times_2 U_2^\dagger \dots \times_d U_d^\dagger) \times_1 U_1 \times_2 U_2 \dots \times_d U_d, \label{tuckerExact:eq}
\end{align}
where the last equality is based on \eqnref{sameModeProdRule:eq}. Formula \eqnref{tuckerExact:eq} implies that if we chose the matrices $U_i$ as the factor matrices and take 
\begin{equation}
\label{TuckerCore:eq}
\mathcal{C} := \mathcal{A} \times_1 U_1^\dagger \times_2 U_2^\dagger \dots \times_d U_d^\dagger
\end{equation}
as an $r_1 \times r_2 \dots \times r_d$ core tensor, then we have the following Tucker decomposition \[
\mathcal{A} = \mathcal{C} \times_1 U_1 \times_2 U_2 \dots \times_d U_d,
\]
where the equality is exact as long as $r_i$ is larger than or equal to the rank of the column space of $A_{(i)}$.

The formulation of HOSVD in \eqnref{hosvd1:eq} is equivalent to
\begin{equation}
\label{multilinear2Kron_supp:eq}
\vec(\mathcal{A}) = (U_d \otimes \dots \otimes U_2 \otimes U_1)\ \vec(\mathcal{C}).
\end{equation}
See \cite[eqn. (12.4.19)]{Golubbookori} for instance.

From \eqnref{multilinear2Kron_supp:eq} and \eqnref{TuckerCore:eq} it is clear that the computation of the core tensor $\mathcal{C}$ is equivalent to solving a huge overdetermined linear system of equations of the form
\[
(U_d \otimes \dots \otimes U_2 \otimes U_1) c = a
\]
where $c:= \vec(\mathcal{C})$ is a vector of size $(r_1 r_2\dots r_d) \times 1$, $a:= \vec(\mathcal{A})$ is a vector of size $N \times 1$ with $N:= n_1 n_2 \dots n_d$ and the coefficient matrix containing Kronecker products is of size $N \times (r_1 r_2\dots r_d)$.

In the case of the deterministic HOSVD, the matrices $U_i$ are chosen to be the left singular vectors of $A_{(i)}$ and the computation of the core relies on the orthogonality of the columns of every $U_i$; see Algorithm \ref{HOSVD:alg}. In the case of the STHOSVD $U_i$ is chosen to be the left singular vectors of the previously-truncated tensor $\hat{\mathcal{C}}^{(i)}$. 

\begin{algorithm}[h!] 
\caption{Deterministic HOSVD (De Lathauwer, De Moor and Vandewalle, 2000 \cite{de2000multilinear}) \\ Inputs are $\mathcal{A} \in \mathbb{R}^{n_1 \times n_2 \times \dots \times n_d}$ and truncation rank $(r_1,r_2, \dots, r_d)$.\\ Output is $\mathcal{A} \approx \llbracket \mathcal{C}; U_1 , U_2 ,\dots, U_d \rrbracket$.}\label{HOSVD:alg}
\begin{algorithmic}[1]
\FOR{$i=1,\ldots,d$}
    \STATE Compute thin SVD $A_{(i)} = \begin{bmatrix}
	\hat U_1 & \hat U_2
	\end{bmatrix}
	\begin{bmatrix}
	\Sigma_1 & \\
	 & \Sigma_2
	\end{bmatrix}
	\begin{bmatrix}
	V_1^T \\
	V_2^T
	\end{bmatrix}$ where $\hat U_1 \in \mathbb{R}^{n_i \times r_i}$.
    \STATE Set $U_i := \hat U_1$.
\ENDFOR
\STATE Compute $\mathcal{C} = \mathcal{A} \times_1 U_1^T \times_2 U_2^T \dots \times_d U_d^T$.
\end{algorithmic}
\end{algorithm}

\begin{algorithm}[h!] 
\caption{Deterministic STHOSVD  (Vannieuwenhoven, Vandebril, and Meerbergen, 2012 \cite{Vannieuwenhoven12})\\ Inputs are $\mathcal{A} \in \mathbb{R}^{n_1 \times n_2 \times \dots \times n_d}$, truncation rank $(r_1,r_2, \dots, r_d)$, and processing order $\bf{p}$ of the modes (a permutation of $[1, 2, \dots , d]$).\\ Output is $\mathcal{A} \approx \llbracket \mathcal{C}; U_1 , U_2 ,\dots, U_d \rrbracket$.}\label{STHOSVD:alg}
\begin{algorithmic}[1]
\STATE Set $\mathcal{C}:= \mathcal{A}$.
\FOR{$i=p_1,\ldots,p_d$}
    \STATE Compute thin SVD $C_{(i)} = \begin{bmatrix}
	\hat U_1 & \hat U_2
	\end{bmatrix}
	\begin{bmatrix}
	\Sigma_1 & \\
	 & \Sigma_2
	\end{bmatrix}
	\begin{bmatrix}
	V_1^T \\
	V_2^T
	\end{bmatrix}$ where $\hat U_1 \in \mathbb{R}^{n_i \times r_i}$.
    \STATE Set $U_i := \hat U_1$.
    \STATE Compute $C_{(i)} = \Sigma_1 V_1^T$.
\ENDFOR
\end{algorithmic}
\end{algorithm}

\section{Tucker to HOSVD conversion}
As mentioned in section~\ref{RHOSVDSMS:subsec}, converting the Tucker decomposition computed via RTSMS to the HOSVD format---where factor matrices possess orthonormal columns and the core tensor is all-orthogonal---is a straightforward process~\cite{de2000multilinear}, based on orthonormalizing the factor matrices with a QR factorization, merging the $R$ factors in the core tensor, and recompressing the updated core tensor for further rank truncation. 
Algorithm~\ref{Tucker2HOSVD:alg} is a standard deterministic approach to accomplishing this conversion.

\begin{algorithm}
\caption{\texttt{Tucker2HOSVD} (with or without thresholding) 
\\ Inputs are Tucker decomposition $\mathcal{A} \approx \llbracket \mathcal{C}; F_1 , F_2 ,\dots, F_d \rrbracket$ whose multilinear rank is $\hat {\bf r}$ (and a tolerance ${\rm tol}$ if thresholding).\\ 
Output is HOSVD $\mathcal{A} \approx  \llbracket \mathcal{\check C}; U_1 , U_2 ,\dots, U_d \rrbracket$ whose multilinear rank is either $\hat {\bf r}$ in case of no thresholding, or ${\bf l}$ in the case of thresholding with $\ell_i \leq \hat r_i$.}\label{Tucker2HOSVD:alg}
\begin{algorithmic}[1]
\FOR{$i=1,2, \ldots,d$}
    \STATE Compute thin QR factorizations $[Q_i, R_i] = qr(F_i)$.
\ENDFOR
\STATE Update $\mathcal{C} := \mathcal{C} \times_1 R_1 \times_2 R_2 \dots \times_d R_d$.
\STATE Apply deterministic STHOSVD to $\mathcal{C}$ and compute $\llbracket \mathcal{\check C}; \check U_1, \check U_2 ,\dots, \check U_d \rrbracket 
\approx \mathcal{C}$ (and mode-$i$ higher order singular values $\sigma^{(i)}$ if thresholding). \hfill \COMMENT{Both $\mathcal{C}$ and $\mathcal{\check C}$ are of size $\hat r_1 \times \hat r_2 \dots \times \hat r_d$.}
\IF{thresholding}
\FOR{$i=1,2, \ldots,d$}  \STATE Find smallest $\ell_i$ such that $\sigma_{\ell_i+1}^{(i)} < {\rm tol}\ \sigma_{1}^{(i)}$.
\ENDFOR
\ELSE
\STATE Set $l_i := \hat r_i$ for $i=1,2,\dots, d$.
\ENDIF
\FOR{$i=1,2, \ldots,d$}  \STATE Compute $U_i =  Q_i\ \check U_i(:, 1:\ell_i)$. \hfill \COMMENT{$Q_i$ is of size $n_i \times \hat r_i$, $\check U_i$ is $\hat r_i \times \hat r_i$, and $U_i$ is $n_i \times \ell_i$.} 
\ENDFOR
\STATE Replace $\mathcal{\check C}$ with $\mathcal{\check C}(1:\ell_1, 1:\ell_2, \dots, 1: \ell_d)$.
\end{algorithmic}
\end{algorithm}

\section{Randomized GN}
Algorithm~\ref{R-GN-ST-Tucker:alg} is a higher-order generalized Nystr\"om for computing a Tucker decomposition in a sequentially truncated manner. It relies on the generalized Nystr\"om framework for randomized low-rank approximation of unfolding matrices as outlined in Algorithm~\ref{RandGN:alg}. It is not the same as the multilinear Nystr\"{o}m algorithm in~\cite{Bucci23}.

\begin{algorithm}[!h] 
\caption{\texttt{GN} Generalized Nystr\"om \cite[Alg. 2.1]{Nakatsukasa20fast}\label{GN:alg}. \\ Inputs are matrix $A \in \mathbb{R}^{m \times n}$, and Gaussian random matrices $\Omega \in \mathbb{R}^{n \times r}$ and $\tilde \Omega \in \mathbb{R}^{m \times \hat r}$.\\ Outputs are $\hat U \in \mathbb{R}^{m \times r}$ and $\hat V \in \mathbb{R}^{n \times r}$ such that $A \approx \hat U \hat V^T$.}\label{RandGN:alg}
\begin{algorithmic}[1]
\STATE Compute $\texttt{AX} := A \Omega$. \STATE Compute $\texttt{YA} := \tilde \Omega^T A$.
\STATE Compute $\texttt{YAX} := \texttt{YA}\ X$.
\STATE Compute thin QR decomposition $\texttt{YAX}  = Q R$. \STATE Compute $\hat U$ by solving $ \hat U R = \texttt{AX}$.
\STATE Compute $\hat V := Q^T\ \texttt{YA}$.
\end{algorithmic}
\end{algorithm}

\begin{algorithm}[!h] 
\caption{\texttt{R-GN-ST-Tucker} \\ Inputs are $\mathcal{A} \in \mathbb{R}^{n_1 \times n_2 \times \dots \times n_d}$, target multilinear rank $(r_1, r_2, \dots, r_d)$, and processing order $\bf{p}$ of the modes.\\ Output is $\mathcal{A} \approx   \llbracket \mathcal{C}; U_1 , U_2 ,\dots, U_d \rrbracket$.}\label{R-GN-ST-Tucker:alg}
\begin{algorithmic}[1]
\STATE Set $\mathcal{C}:= \mathcal{A}$.
\FOR{$i=p_1,\ldots,p_d$}
    \STATE Draw two standard random Gaussian matrices $\Omega_i$ and $\tilde \Omega_i$ of size $z_{i} \times r_{i}$ and $n_{i} \times \hat r_{i}$, respectively, where $\hat r_{i}:= r_{i} +p$ with $p := [ r_i/2 ]$.
    \STATE Compute $[\hat U, \hat V] = \texttt{GN}(C_{(i)}, \Omega_i, \tilde \Omega_i)$ using Alg.~\ref{GN:alg}.
    \STATE Set $U_i := \hat U$. 
    \STATE Update $C_{(i)} = \hat V^T$. \hfill       \COMMENT{Overwriting $C_{(i)}$ overwrites $\mathcal{C}$.}
\ENDFOR
\end{algorithmic}
\end{algorithm}

\end{document}